\documentclass[11pt]{amsart}

\usepackage{amssymb,graphicx}

\input xy
\xyoption{all}

\newtheorem{proposition}{Proposition}
\newtheorem{theorem}[proposition]{Theorem}
\newtheorem{lemma}[proposition]{Lemma}
\newtheorem{corollary}[proposition]{Corollary}
\newtheorem{conjecture}[proposition]{Conjecture}

\theoremstyle{definition}
\newtheorem{definition}[proposition]{Definition}

\theoremstyle{remark}

\numberwithin{proposition}{section}

\def\lk{\operatorname{lk}}
\def\ring{\Z[\lambda^{\pm 1},\mu^{\pm 1}]}
\def\Diag{\Delta}
\def\k{\mathbf{k}}
\def\vareps{\varepsilon}
\def\eps{\epsilon}
\def\X{X}

\def\Z{\mathbb{Z}}
\def\R{\mathbb{R}}
\def\Q{\mathbb{Q}}
\def\A{\mathcal{A}}

\def\I{\mathcal{I}}

\def\M{\mathcal{M}}

\def\d{\partial}

\def\Aut{\operatorname{Aut}}

\def\hom#1{\phi_{#1}}
\def\Phil#1{\Phi^L_{#1}}
\def\Phir#1{\Phi^R_{#1}}

\def\Bur#1{{\operatorname{Bur}_{#1}}}

\def\Aug{\operatorname{Aug}}

\def\ext{\operatorname{ext}}

\def\diag{\operatorname{diag}}
\def\lin{\operatorname{lin}}

\def\WB{\textit{WB}}

\def\C{\mathbb{C}}

\def\int{\operatorname{int}}
\def\P{\mathcal{P}}
\def\secdivide{}
\def\subsecdivide{}

\theoremstyle{plain}

\begin{document}

\title{Framed Knot Contact Homology}
\author[L. Ng]{Lenhard Ng}
\address{Department of Mathematics, Duke
University, Durham, NC 27708} \email{ng@math.duke.edu}
\urladdr{http://alum.mit.edu/www/ng}

\begin{abstract}
We extend knot contact homology to a theory over the ring
$\mathbb{Z}[\lambda^{\pm 1},\mu^{\pm 1}]$, with the invariant given
topologically and combinatorially. The improved invariant, which is
defined for framed knots in $S^3$ and can be generalized to knots in
arbitrary manifolds, distinguishes the unknot and can distinguish
mutants. It contains the Alexander polynomial and naturally produces
a two-variable polynomial knot invariant which is related to the
$A$-polynomial.
\end{abstract}

\maketitle

%*********************************************************************
%*********************************************************************
\section{Introduction}
\label{sec:intro}

This may be viewed as the third in a series of papers on knot
contact homology, following \cite{Ng1,Ng2}. In this paper, we extend
the knot invariants of \cite{Ng1,Ng2}, which were defined over the
base ring $\Z$, to invariants over the larger ring $\ring$. The new
invariants, defined both for knots in $S^3$ or $\R^3$ and for more
general knots, contain a large amount of information not contained
in the original formulation of knot contact homology, which
corresponds to specifying $\lambda=\mu=1$. We will recapitulate
definitions from the previous papers where necessary, so that the
results from this paper can hopefully be understood independently of
\cite{Ng1,Ng2}, although some proofs rely heavily on the previous
papers.

The motivation for the invariant given by knot contact homology
comes from symplectic geometry and in particular the Symplectic
Field Theory of Eliashberg, Givental, and Hofer \cite{EGH}. Our
approach in this paper and its predecessors is to view the invariant
topologically and combinatorially, without using the language of
symplectic topology. Work is currently in progress to show that the
invariant defined here does actually give the ``Legendrian contact
homology'' of a certain canonically defined object. For more
details, see \cite[\S 3]{Ng1} or the end of this section, where we
place the current work in context within the general subject of
holomorphic curves in symplectic manifolds.

We introduce the invariant in Section~\ref{sec:defs}. For a knot in
$\R^3$, the form for the invariant which carries the most
information is a differential graded algebra over $\ring$, which we
call the \textit{framed knot DGA}, modulo an equivalence relation
known in contact geometry as stable tame isomorphism, which is a
special case of quasi-isomorphism. We can define this algebra in
terms of a braid whose closure is the desired knot; this is the
version of the invariant which is directly derived from
considerations in contact geometry. There is a slightly more natural
definition for the framed knot DGA in terms of a diagram of the
knot, but a satisfying topological interpretation for the full
invariant has yet to be discovered.

We refer to the homology of the framed knot DGA, which is also a
knot invariant, as \textit{framed knot contact homology}. The most
tractable piece of information which can be extracted from the
framed knot DGA is the degree $0$ piece of this homology, which we
call the \textit{cord algebra} of the knot because it affords a
natural topological reformulation in terms of cords, in the spirit
of \cite{Ng2}. The cord algebra can easily be extended to knots in
arbitrary manifolds using homotopy groups, specifically the
peripheral information attached to the knot group. It seems likely
that this general cord algebra measures the degree $0$ part of the
corresponding Legendrian contact homology, as for knots in $S^3$,
but a proof would need to expand our current technology for
calculating contact homology.

The framed knot DGA contains a fair amount of ``classical''
information about the knot. In particular, one can deduce the
Alexander polynomial from a certain canonical linearization of the
algebra. More interestingly, one can use the framed knot DGA, or
indeed the cord algebra derived from it, to define a two-variable
polynomial knot invariant which we call the \textit{augmentation
polynomial}. The augmentation polynomial, in turn, has a factor
which is the well-studied $A$-polynomial introduced in \cite{CCGLS}.
In particular, a result of Dunfield and Garoufalidis \cite{DG} about
the nontriviality of the $A$-polynomial implies that the cord
algebra (and hence framed knot contact homology) distinguishes the
unknot from all other knots. There is a close relationship between
knot contact homology and $SL_2\C$ representations of the knot
complement, as evidenced by the link to the $A$-polynomial, but
there does not seem to be an interpretation for the framed knot DGA
purely in terms of character varieties.

On a related topic, the cord algebra is a reasonably effective tool
to distinguish knots. Certainly it can tell apart knots which have
different Alexander polynomials, and it can also distinguish mirrors
(e.g., left handed and right handed trefoils). It can even
distinguish knots which have the same $A$-polynomial, such as the
Kinoshita--Terasaka knot and its Conway mutant. The tool used to
produce such results is a collection of numerical invariants derived
from the cord algebra, called \textit{augmentation numbers}, which
can be readily calculated by computer. It is even remotely possible
that the cord algebra could be a complete knot invariant.

Here is an outline of the paper. Section~\ref{sec:defs} contains the
various definitions of the knot invariant, and the nontrivial
portions of the proofs of their invariance and equivalence are given
in Section~\ref{sec:proofs}. In Section~\ref{sec:properties}, we
discuss properties of the invariant, including the aforementioned
relation with the Alexander polynomial. Section~\ref{sec:aug}
examines some information which can be extracted from the invariant,
including augmentation numbers and the augmentation polynomial, and
discusses the relation with the $A$-polynomial.
Section~\ref{sec:extensions} outlines extensions of the invariant to
various other contexts, including spatial graphs and
higher-dimensional knots.

We conclude this section with some remarks about the contact
geometry of the framed knot DGA. As discussed in \cite{Ng1}, the
knot DGA of a knot in $\R^3$ measures the Legendrian contact
homology of a naturally defined Legendrian torus (the unit conormal
to the knot) in the contact manifold $ST^*\R^3$. In its most general
form (see \cite{ENS}), Legendrian contact homology is defined over
the group ring of the first homology group of the Legendrian
submanifold: briefly, the group ring coefficients arise from
considering the boundaries of the holomorphic disks used in
Legendrian contact homology, which lie in the Legendrian submanifold
and can be closed to closed curves by adding capping paths joining
endpoints of Reeb chords.

In our case, the Legendrian torus is essentially the boundary of a
tubular neighborhood of the knot, and the group ring $\Z H_1(T^2)$
can be identified with $\ring$ once we choose a framing for the
knot. The resulting DGA, whose homology is the Legendrian contact
homology of the Legendrian torus, is precisely the framed knot DGA
defined above; this result is the focus of an upcoming paper of
Ekholm, Etnyre, Sullivan, and the author.

Since its introduction by Eliashberg and Hofer \cite{Eli},
Legendrian contact homology has been studied in many papers,
including \cite{Che,EES,ENS,Ng}. Much of this work has focused on
the lowest dimensional case, that of Legendrian knots in standard
contact $\R^3$ (or other contact three-manifolds). The present
manuscript, along with its predecessors \cite{Ng1,Ng2}, comprises
one of the first instances of a reasonably nontrivial and involved
computation of Legendrian contact homology in higher dimensions; see
also \cite{EES3}. (This statement tacitly assumes the pending
result, mentioned above, that the combinatorial version of knot
contact homology presented here agrees with the geometric version
from contact topology.)

Beyond its implications for knot theory, the Legendrian contact
homology in this case is interesting because it has several features
not previously observed in the subject. In particular, the use of
group ring coefficients in our setting is crucial for applications,
whereas the author currently knows of no applications of group ring
coefficients for, say, the theory of one-dimensional Legendrian
knots.
%We have seen
%in \cite{Ng1} that the knot DGA over $\Z_2$ carries augmentations
%which give rise to distinct Poincar\'e polynomials; to the author's
%knowledge, this had never been seen before, and in fact it may never
%be possible for one-dimensional Legendrian knots. In the present
%paper,

%\begin{figure}
%\[
%\xymatrix{ \textrm{Homotopy invariant} \ar[d] & \textrm{Framed knot
%DGA} \ar@/_1pc/[dl] \ar[d] \ar[r]
%& \textrm{Knot DGA}/\Z \ar[d] \ar@/^4pc/[dd] \\
%HC_0=\textrm{cord algebra} \ar[dr] \ar[d] \ar@/^2pc/[rr] \ar[drr] & HC_*^{\lin} \ar@/^5pc/[dd]
%& \widetilde{HC}_0,\widetilde{HC}_0^{\ab}/\Z \\
%\textrm{Augmentation variety} \ar[d] \ar@{-->}[dr] & \textrm{Augmentation numbers}
%& HC_0,HC_0^{\ab}/\Z \ar@{-->}[d] \\
%\textrm{Slope polynomial ($\sigma$?)} & \Delta_K(t) & \textrm{Character variety}
%}
%\]
%\caption{
%The invariants discussed in this paper, and their interrelationships.
%Solid arrows are direct derivations; dashed arrows are
%implicit relations.
%}
%\label{fig:summary}
%\end{figure}

\subsecdivide
%*********************************************************************
\subsection*{Acknowledgments}

I would like to thank Tobias Ekholm, Yasha Eliashberg, John Etnyre,
Tom Mrowka, Colin Rourke, Mike Sullivan, and Ravi Vakil for useful
discussions. This work was supported by an American Institute of
Mathematics Five-Year Fellowship and NSF grant FRG-0244663.

\secdivide
%*********************************************************************
%*********************************************************************
\section{The Invariant}
\label{sec:defs}

We will give four interpretations of the invariant. The first, in
terms of homotopy groups, generalizes the version of the cord ring
given in the appendix to \cite{Ng2}. The second, in terms of cords,
generalizes the original definition of the cord ring in \cite{Ng2}.
The third, in terms of a braid representation of the knot, and the
fourth, in terms of a knot diagram, extend the definition of knot
contact homology and the knot DGA from \cite{Ng1}.

Throughout this section, we work with noncommutative algebras, but
for the purposes of our applications, one could just as well
abelianize and work in the commutative category. (Note that for the
full DGA, abelianizing forces us to work over a ring containing
$\mathbb{Q}$ \cite{ENS}.) This is in contrast to the situation of
Legendrian knots in $\R^3$, where it is sometimes advantageous to
exploit noncommutativity \cite{Ng}. It is possible that one could
use the noncommutative structure for the cord algebra or framed knot
DGA to distinguish a knot from its inverse; see
Section~\ref{ssec:symmetries}.

\subsecdivide
%*********************************************************************
\subsection{Homotopy interpretation}
\label{ssec:htpydef}

This version of the invariant has the advantage of being defined in
the most general setup.

Let $K\subset M$ be a connected codimension $2$ submanifold of any
smooth manifold, equipped with an orientation of the normal bundle
to $K$ in $M$. Denote by $\nu(K)$ a tubular neighborhood of $K$ in
$M$, and choose a point $x_0$ on the boundary $\d(\nu(K))$. The
orientation of the normal bundle determines an element $m \in
\pi_1(\d\nu K,x_0)$ given by the oriented fiber of the $S^1$ bundle
$\d\nu K$ over $K$. The inclusion $\iota:\thinspace\d\nu K
\hookrightarrow M\setminus K$ induces a map $\iota_*:\thinspace
\pi_1(\d\nu K,x_0) \rightarrow \pi_1(M\setminus K,x_0)$.

Denote by $R=\Z H_1(\d\nu K)$ the group ring of the homology group
$H_1(\d\nu K)$, with coefficients in $\Z$; the usual abelianization
map from $\pi_1$ to $H_1$ yields a map $p:\thinspace \pi_1(\d\nu
K,x_0) \rightarrow R$. Let $\A$ denote the tensor algebra over $R$
freely generated by the elements of the group $\pi_1(M\setminus
K,x_0)$, where we view $\pi_1(M\setminus K,x_0)$ as a set. We write
the image of $\gamma\in\pi_1(M\setminus K,x_0)$ in $\A$ as
$[\gamma]$; then $\A$ is generated as an $R$-module by words of the
form $[\gamma_1]\cdots[\gamma_k]$, $k\geq 0$. (Note that
$[\gamma_1\gamma_2] \neq [\gamma_1][\gamma_2]$ in $\A$.)

\begin{definition}
Let $\I\subset\A$ be the subalgebra generated by the following elements of $\A$:
\begin{enumerate}
\item $[e]-1-p(m)$, where $e$ is the identity in $\pi_1(M\setminus K)$ and $1$ is the
unit in $\A$;
\item $[\gamma\iota_*(\gamma')]-p(\gamma') [\gamma]$ and
$[\iota_*(\gamma')\gamma]-p(\gamma') [\gamma]$, for any $\gamma'\in\pi_1(\d\nu K)$ and
$\gamma\in\pi_1(M\setminus K)$;
\item $[\gamma_1 \gamma_2] + [\gamma_1 \iota_*(m) \gamma_2] - [\gamma_1] \cdot [\gamma_2]$,
for any $\gamma_1,\gamma_2\in\pi_1(M\setminus K)$.
\end{enumerate}
The \textit{cord algebra} of $K\subset M$ is the algebra $\A/\I$, over the ring
$R$.
\label{def:generalcordalg}
\end{definition}

This algebra, modulo isomorphism fixing the base ring $R$, is
evidently an isotopy invariant of $K$. We note that the relation
$[e]=1+p(m)$ in $\A/\I$ is chosen to give consistency between the
other two relations; the third relation gives
\[
[\gamma_1][e] = [\gamma_1 e] + [\gamma_1 \iota_*(m) e] =
[\gamma_1] + [\gamma_1 \iota_*(m)]
\]
for $\gamma_1 \in \pi_1(M\setminus K)$,
while the second relation gives $[\gamma_1] + [\gamma_1 \iota_*(m)] = (1+p(m)) [\gamma_1]$.

In our case of interest, $K$ is an oriented knot in $M = S^3$ (or
$\R^3$), $\d\nu K$ is topologically a $2$-torus, and $H_1(\d\nu K) =
\Z^2$. A framing of $K$ gives a set $\{\lambda,\mu\}$ of generators
of $H_1(\d\nu K)$, where $\lambda$ is the longitude of $K$ and $\mu$
is the meridian. Given a framing, we can hence identify $R = \Z
H_1(\d\nu K)$ with the ring $\Z[\lambda^{\pm 1},\mu^{\pm 1}]$.

\begin{definition}
Let $K\subset S^3$ be a framed knot, and let $l,m$ denote the
homotopy classes of the longitude and meridian of $K$ in
$\pi_1(S^3\setminus K)$. The \textit{framed cord algebra} of $K$ is
the tensor algebra over $\Z[\lambda^{\pm 1},\mu^{\pm 1}]$ freely
generated by $\pi_1(S^3\setminus K)$, modulo the relations
\begin{enumerate}
\item \label{htpy1} $[e] = 1+\mu$;
\item \label{htpy2} $[\gamma m] = [m\gamma] = \mu [\gamma]$ and $[\gamma l] = [l\gamma] = \lambda [\gamma]$,
for $\gamma\in\pi_1(S^3\setminus K)$;
\item \label{htpy3} $[\gamma_1 \gamma_2] + [\gamma_1 m \gamma_2] = [\gamma_1][\gamma_2]$,
for $\gamma_1,\gamma_2\in\pi_1(S^3\setminus K)$.
% changed this!
\end{enumerate}
\label{def:knotcordalg}
\end{definition}

\noindent It is clear that, for framed knots, this definition agrees
with Definition~\ref{def:generalcordalg} above.

A few remarks are in order. First, the condition $[\gamma m] =
[m\gamma] = \mu [\gamma]$ in (\ref{htpy2}) of
Definition~\ref{def:knotcordalg} is actually unnecessary since it
follows from (\ref{htpy1}) and (\ref{htpy3}) with
$(\gamma_1,\gamma_2) = (\gamma,e)$ or $(e,\gamma)$.

Second, changing the framing of the knot by $k\in\Z$ has the effect
of replacing $\lambda$ by $\lambda \mu^k$ in the cord algebra. If we
change $K$ to the mirror of $K$, with the corresponding framing,
then the cord algebra changes by replacing $\mu$ by $\mu^{-1}$; if
we reverse the orientation of $K$, then the cord algebra changes by
replacing $\lambda$ by $\lambda^{-1}$ and $\mu$ by $\mu^{-1}$.

Finally, as an illustration of Definition~\ref{def:knotcordalg}, we
compute the framed cord algebra for the $0$-framed unknot in $S^3$.
Here $\pi_1(S^3\setminus K) \cong \Z$ is generated by the meridian
$m$, and the relations in Definition~\ref{def:knotcordalg} imply
that $[e]=1+\mu$, $[m]=\mu[e]$, and $[e]=\lambda[e]$, and so the
framed cord algebra is
\[
\ring/((\lambda-1)(\mu+1)).
\]

\subsecdivide
%*********************************************************************
\subsection{Cord interpretation}
\label{ssec:corddef}

The homotopy definition of the cord algebra given above is perhaps
not the most useful formulation from a computational standpoint.
Here we give another topological definition of the cord algebra, in
terms of the paths which we termed ``cords'' in \cite{Ng2}.

\begin{definition} \label{def:cord}
Let $K\subset S^3$ be an oriented knot, and let $*$ be a fixed point
on $K$. A \textit{cord} of $(K,*)$ is any continuous path
$\gamma:\thinspace [0,1]\rightarrow S^3$ with $\gamma^{-1}(K) =
\{0,1\}$ and $\gamma^{-1}(\{*\}) = \emptyset$. Let $\A_K$ denote the
tensor algebra over $\Z[\lambda^{\pm 1},\mu^{\pm 1}]$ freely
generated by the set of homotopy classes of cords of $(K,*)$.
\end{definition}

Note that this definition differs slightly from the definition of
cords in \cite{Ng2}. As in \cite{Ng2}, we distinguish between a knot
and its cords in diagrams by drawing the knot more thickly.

In $\A_K$, we define skein relations as follows:
\begin{equation}
\label{eq:skein1new}
\raisebox{-0.17in}{\includegraphics[width=0.4in]{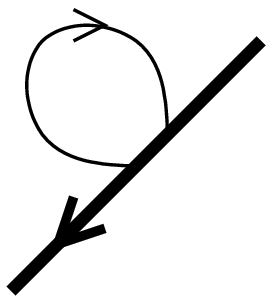}} =
1+\mu;
\end{equation}
\begin{equation}
\label{eq:skein2new}
\raisebox{-0.17in}{\includegraphics[width=0.4in]{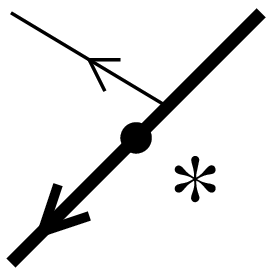}} =
\lambda
\raisebox{-0.17in}{\includegraphics[width=0.4in]{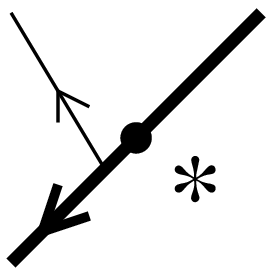}}
\hspace{0.4in} \textrm{and} \hspace{0.4in}
\raisebox{-0.17in}{\includegraphics[width=0.4in]{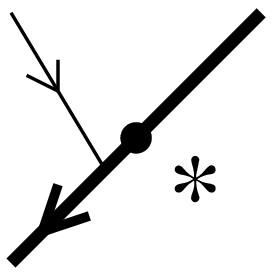}} =
\lambda
\raisebox{-0.17in}{\includegraphics[width=0.4in]{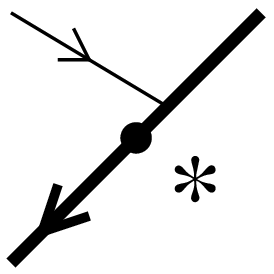}};
\end{equation}
\begin{equation}
\label{eq:skein3new}
\raisebox{-0.17in}{\includegraphics[width=0.4in]{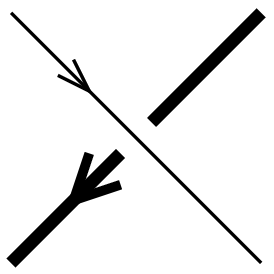}} + \mu
\raisebox{-0.17in}{\includegraphics[width=0.4in]{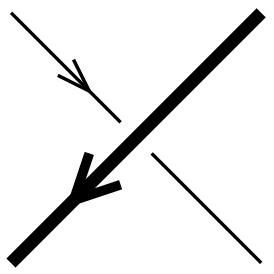}} =
\raisebox{-0.17in}{\includegraphics[width=0.4in]{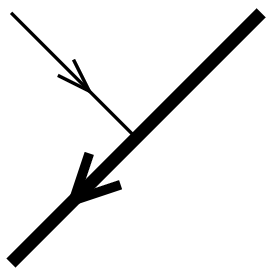}} \cdot
\raisebox{-0.17in}{\includegraphics[width=0.4in]{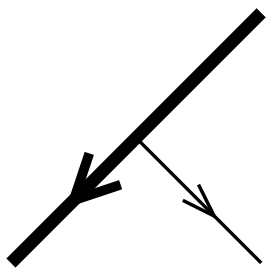}}.
\end{equation}
These diagrams are understood to depict some neighborhood in $S^3$
outside of which the diagrams agree. To clarify the diagrams, if one
pushes either of the cords on the left hand side of
(\ref{eq:skein3new}) to intersect the knot, then the cord splits
into the two cords on the right hand side of (\ref{eq:skein3new});
also, the cord in (\ref{eq:skein1new}) is any contractible cord.
Note that the skein relations are considered as diagrams in space
rather than in the plane. For example, (\ref{eq:skein3new}) is
equivalent to
\[
\raisebox{-0.17in}{\includegraphics[width=0.4in]{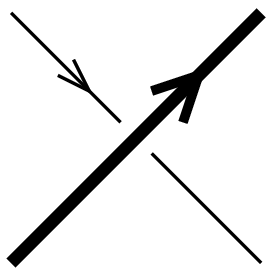}} +
\mu
\raisebox{-0.17in}{\includegraphics[width=0.4in]{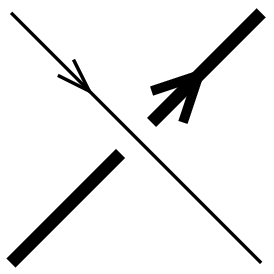}} =
\raisebox{-0.17in}{\includegraphics[width=0.4in]{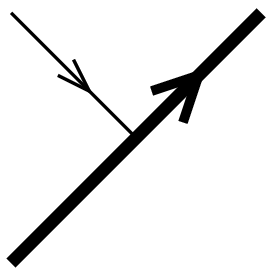}}
\cdot
\raisebox{-0.17in}{\includegraphics[width=0.4in]{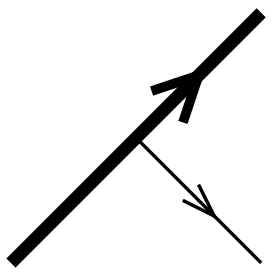}},
\]
even though the two relations are not the same when interpreted as plane diagrams.

\begin{definition}
The \textit{cord algebra} of $(K,*)$ is $\A_K$ modulo the relations
(\ref{eq:skein1new}), (\ref{eq:skein2new}), and
(\ref{eq:skein3new}). \label{def:cordcordalg}
\end{definition}

Note that the definition of the cord algebra of $(K,*)$, if we set
$\lambda=\mu=1$, is the same as the definition of the cord ring of
$K$ from \cite{Ng2}, after we change variables by replacing every
cord by $-1$ times itself.

\begin{proposition}
The cord algebra of $(K,*)$ is independent, up to isomorphism, of
the choice of the point $*$, and is isomorphic to the framed cord
algebra of Definition~\ref{def:knotcordalg} when $K\subset S^3$ is
given the $0$ framing.
\end{proposition}

\begin{proof}
To see the isomorphism between the algebras specified by
Definitions~\ref{def:knotcordalg} and \ref{def:cordcordalg}, connect
the base point $x_0 \in \d\nu K$ from the homotopy interpretation to
a nearby point $x\in K$ (with $x\neq *$) via some path. By
conjugating by this path, we can associate, to any loop $\gamma$
representing an element of $\pi_1(S^3\setminus K,x_0)$, a cord
$\tilde{\gamma}$ beginning and ending at $x$, and homotopic loops
are mapped to homotopic cords. The map from the framed cord algebra
with framing $0$ to the cord algebra of
Definition~\ref{def:cordcordalg} sends $[\gamma]$ to
$\mu^{\lk(\gamma,K)} [\tilde{\gamma}]$, where $\lk$ is the linking
number. It is easily verified that this map is an isomorphism. Note
that the $0$ framing is necessary so that a cord following the
longitude of $K$ is identified with $\lambda$ times the contractible
cord.
\end{proof}

It follows readily from Definition~\ref{def:cordcordalg} that the
cord algebra of a knot $K$ is finitely generated and finitely
presented. Any cord, in fact, can be expressed in terms of the
so-called minimal binormal chords of the knot, which are the local
length minimizers among the line segments beginning and ending on
the knot (imagine a general cord as a rubber band and pull it
tight). For a description of how we can calculate the cord algebra
from a diagram of the knot, see Section~\ref{ssec:cordDGA}.

\subsecdivide
%*********************************************************************
\subsection{Differential graded algebra interpretation I}
\label{ssec:DGAdef}

This version of the invariant is the original one, derived from
holomorphic curve theory and contact homology. The invariant we
define here is actually ``better'' than the cord algebra defined
above; it is a differential graded algebra whose homology in degree
$0$ is the cord algebra. It is an open question to develop a
topological interpretation for the full differential graded algebra,
using cords or homotopy groups as in Sections \ref{ssec:htpydef} and
\ref{ssec:corddef}.

We recall some definitions from \cite{Ng1}, slightly modified for
our purposes. Let $B_n$ be the braid group on $n$ strands, and let
$\A_n$ denote the tensor algebra over $\ring$ freely generated by
the $n(n-1)$ generators $a_{ij}$, with $1\leq i,j\leq n$ and $i\neq
j$. There is a representation $\phi:\thinspace
B_n\rightarrow\Aut\A_n$ defined on the usual generators
$\sigma_1,\ldots, \sigma_{n-1}$ of $B_n$ as follows:
\[
\hom{\sigma_k}:\thinspace \left\{
\begin{array}{ccll}
a_{ki} & \mapsto & -a_{k+1,i} - a_{k+1,k}a_{ki} & i\neq k,k+1 \\
a_{ik} & \mapsto & -a_{i,k+1} - a_{ik}a_{k,k+1} & i\neq k,k+1 \\
a_{k+1,i} & \mapsto & a_{ki} & i\neq k,k+1 \\
a_{i,k+1} & \mapsto & a_{ik} & i \neq k,k+1 \\
a_{k,k+1} & \mapsto & a_{k+1,k} & \\
a_{k+1,k} & \mapsto & a_{k,k+1} & \\
a_{ij} & \mapsto & a_{ij} & i,j \neq k,k+1.
\end{array}
\right.
\]
Write $\phi_B$ as the image of a braid $B\in B_n$ under $\phi$.
There is a similar representation $\phi^{\ext}:\thinspace
B_n\rightarrow\Aut\A_{n+1}$ given by the composition of the
inclusion $B_n \hookrightarrow B_{n+1}$ with the map $\phi$ on
$B_{n+1}$, where the inclusion adds a strand we label $*$ to the
$n$-strand braid. For $B\in B_n$, we can then define $n\times n$
matrices $\Phil{B},\Phir{B}$ with coefficients in $\A_n$ by the
defining relations
\[
\hom{B}^{\ext}(a_{i*}) = \sum_{j=1}^n (\Phil{B})_{ij} a_{j*}
\hspace{0.25in} \textrm{and} \hspace{0.25in}
 \hom{B}^{\ext}(a_{*j}) =
\sum_{i=1}^n a_{*i} (\Phir{B})_{ij}.
\]

Fix a braid $B\in B_n$. Let $\A$ be the graded tensor algebra over
$\ring$ generated by: $a_{ij}$, $1\leq i,j\leq n$, $i\neq j$, of
degree $0$; $b_{ij}$ and $c_{ij}$, $1\leq i,j\leq n$, of degree $1$;
and $d_{ij}$, $1\leq i,j\leq n$, and $e_i$, $1\leq i\leq n$, of
degree $2$. (Note that $\A$ contains $\A_n$ as a subalgebra.)
Assemble the $a,b,c,d$ generators of $\A$ into $n\times n$ matrices
as follows: $B=(b_{ij})$, $C=(c_{ij})$, $D=(d_{ij})$, and
$A=(a_{ij}')$, where
\[
a_{ij}' = \begin{cases}
\mu a_{ij}, & i<j \\
a_{ij}, & i>j \\
-1-\mu, & i=j.
\end{cases}
\]
(The matrix $B$ is only tangentially related to the braid $B$; the
distinction should be clear from context.) Also define $\Lambda$ to
be the $n\times n$ diagonal matrix $\diag(\lambda,1,\ldots,1)$.

\begin{definition}
Let $B\in B_n$, and let $\A$ be the algebra given above. Define a
differential $\d$ on the generators of $\A$ by
\begin{eqnarray*}
\d A &=& 0 \\
\d B &=& (1 - \Lambda \cdot \Phil{B}) \cdot A \\
\d C &=& A \cdot (1 - \Phir{B}\cdot\Lambda^{-1}) \\
\d D &=& B \cdot (1 - \Phir{B}\cdot\Lambda^{-1}) - (1 - \Lambda\cdot\Phil{B}) \cdot C \\
\d e_i &=& (B + \Lambda\cdot\Phil{B} \cdot C)_{ii},
\end{eqnarray*}
where $\cdot$ denotes matrix multiplication, and extend $\d$ to $\A$
via linearity over $\ring$ and the (signed) Leibniz rule. Then
$(\A,\d)$ is the \textit{framed knot DGA} of $B$.
\label{def:knotDGAbraid}
\end{definition}

Here ``DGA'' is the abbreviation commonly used in the subject for a
(semifree) differential graded algebra. We remark that if we set
$\lambda=\mu=1$, we recover the definition of the knot DGA over $\Z$
from \cite{Ng1}.

There is a standard notion of equivalence between DGAs known as
\textit{stable tame isomorphism}, originally due to Chekanov
\cite{Che}; we now briefly review its definition, in the version
which we need. Note that an important property of stable tame
isomorphism is that it preserves homology (see \cite{ENS}).

Suppose that we have two DGAs $(\A,\d)$ and $(\A',\d')$, where
$\A,\A'$ are tensor algebras over $\ring$ generated by
$a_1,\ldots,a_n$ and $a_1',\ldots,a_n'$, respectively. An algebra
map $\psi:\thinspace \A\rightarrow\A'$ is an \textit{elementary
isomorphism} if it is a graded chain map and, for some $i$,
\begin{align*}
\psi(a_i) &= \alpha a_i'+v  \\
\psi(a_j) &= a_j' \hspace{6ex} \text{ for all } j\neq i,
\end{align*}
where $\alpha$ is a unit in $\ring$ and $v$ is in the subalgebra of
$\A'$ generated by $a_1',\ldots,a_{i-1}',a_{i+1}',\ldots,a_n'$. A
\textit{tame isomorphism} between DGAs is a composition of
elementary isomorphisms. Let $(\mathcal{E}^i,\d^i)$ be the tensor
algebra on two generators $e_1^i,e_2^i$, with $\deg e_1^i-1=\deg
e_2^i=i$ and differential induced by $\d^ie_1^i=e_2^i, \d^i
e_2^i=0$. The degree-$i$ \textit{algebraic stabilization} of a DGA
$(\A,\d)$ is the coproduct of $\A$ with $\mathcal{E}^i$, with the
differential induced from $\d$ and $\d^i$. Finally, two DGAs are
\textit{stable tame isomorphic} if they are tamely isomorphic after
some (possibly trivial, possibly different) number of algebraic
stabilizations of each.

We can now state the invariance result for framed knot DGAs. For any
braid $B\in B_n$, define the \textit{writhe} $w(B)$, as usual, to be
the sum of the exponents of the word in
$\sigma_1,\ldots,\sigma_{n-1}$ comprising $B$.

\begin{theorem}
If braids $B$ and $B'$ have isotopic knot closures, then the framed
knot DGA for $B'$ is stable tame isomorphic to the DGA obtained by
replacing $\lambda$ by $\lambda \mu^{w(B')-w(B)}$ in the framed knot
DGA for $B$.
\label{thm:DGAinv}
\end{theorem}

\noindent We defer the proof of Theorem~\ref{thm:DGAinv} to
Section~\ref{ssec:pfDGAinv}.

\begin{definition}
Let $K\subset S^3$ be a knot. The \textit{$f$-framed knot DGA of
$K$} is the (stable tame isomorphism class of the) framed knot DGA
of any braid $B$ whose closure is $K$ and which has writhe $f$. We
will often refer to the $0$-framed knot DGA, which can be obtained
from the $f$-framed knot DGA by replacing $\lambda$ by
$\lambda\mu^f$, as the \textit{framed knot DGA}.
\label{def:knotDGA}
\end{definition}

%\noindent
%We will mainly be interested in the $0$-framed knot DGA; note that Theorem~\ref{thm:DGAinv} states
%that the $f$-framed knot DGA is precisely the $0$-framed knot DGA with $\lambda$ replaced by
%$\lambda\mu^f$.

\noindent Since mirroring affects the framed knot DGA nontrivially,
Definition~\ref{def:knotDGA} is not completely clear until we
specify what positive and negative braid crossings look like when
the braid is embedded in $\R^3$. We choose the convention that a
positive crossing in a braid (i.e., a generator $\sigma_k$) is one
in which the strands rotate around each other counterclockwise. Thus
the writhe of a braid is \textit{negative} the writhe of the knot
which is its closure; we will henceforth only use the term
``writhe'' vis-\`a-vis braids.

Since stable tame isomorphic DGAs have isomorphic homology, the
homology of the framed knot DGA produces another framed knot
invariant.

\begin{definition}
Let $K$ be a framed knot in $S^3$ with framing $f$, or a knot with
framing $0$ if no framing is specified. The \textit{framed knot
contact homology} of $K$, denoted $HC_*(K)$, is the graded homology
of the framed knot DGA of $K$.
\label{def:framedHC}
\end{definition}

Note that $HC_*(K)$ only exists in dimensions $*\geq 0$. The
following result connects the invariant defined in this section to
the invariants from the two previous sections.

\begin{theorem}
Let $K\subset S^3$ be a knot. Then the degree $0$ framed knot
contact homology $HC_0(K)$, viewed as an algebra over $\ring$, is
isomorphic to the cord algebra of $K$.
\label{thm:DGAcord}
\end{theorem}

\noindent We will indicate a proof of Theorem~\ref{thm:DGAcord} in
Section~\ref{ssec:pfDGAcord}.
%Note the following analogue of
%Proposition 4.2 in \cite{Ng1}: if $B\in B_n$ is a braid with writhe
%$f$ whose closure is a knot $K$, then $HC_0(K)$ is the algebra
%$\A_n$ modulo the subalgebra generated by the entries of the
%matrices $(1 - \Lambda \cdot \Phil{B}) \cdot A$ and $A \cdot (1 -
%\Phir{B}\cdot\Lambda^{-1})$.

\begin{corollary}[{cf.\ \cite[Prop.\ 4.2]{Ng1}}]
If $B\in B_n$ is a braid which has closure $K$, then the cord
algebra of $K$ is the result of replacing $\lambda$ by
$\lambda\mu^{-w(B)}$ in $\A_n/\I$, where $\I$ is the subalgebra
generated by the entries of the matrices $(1 - \Lambda \cdot
\Phil{B}) \cdot A$ and $A \cdot (1 - \Phir{B}\cdot\Lambda^{-1})$.
\end{corollary}

Potentially the framed knot DGA of a knot contains much more
information than the cord algebra. However, it seems difficult to
compute $HC_*$ for $*>0$, and it remains an open question whether
stable tame isomorphism is a stronger relation than
quasi-isomorphism (i.e., having isomorphic homology).

As an example, we compute the cord algebra for the left handed
trefoil, which is the closure of the braid $\sigma_1^3\in B_2$. As
calculated in \cite{Ng1}, we have
%\[
%\Phil{\sigma_1^{-3}} = \left(
%\begin{matrix} a_{21} & -1+a_{21}a_{12} \\
%1-a_{12}a_{21} & 2a_{12}-a_{12}a_{21}a_{12}
%\end{matrix} \right)
%\quad \textrm{and} \quad
%\Phir{\sigma_1^{-3}} = \left(
%\begin{matrix} a_{12} & 1-a_{12}a_{21} \\
%-1+a_{21}a_{12} & 2a_{21}-a_{21}a_{12}a_{21}
%\end{matrix} \right) ,
%\]
\[
\Phil{\sigma_1^3} = \left(
\begin{smallmatrix} 2a_{21}-a_{21}a_{12}a_{21} & 1-a_{21}a_{12} \\
-1+a_{12}a_{21} & a_{12}
\end{smallmatrix} \right)
\quad \textrm{and} \quad
\Phir{\sigma_1^3} = \left(
\begin{smallmatrix} 2a_{12}-a_{12}a_{21}a_{12} & -1+a_{12}a_{21} \\
1-a_{21}a_{12} & a_{21}
\end{smallmatrix} \right) ,
\]
and hence $HC_0$ is $\A_2$ modulo the entries of the matrices
\[
(1-\Lambda\cdot\Phil{\sigma_1^{3}})\cdot A = \left(
\begin{smallmatrix} 
-1-\mu+\lambda(1+2\mu)a_{21}-\lambda\mu a_{21}a_{12}a_{21} &
\lambda+\lambda\mu+\mu a_{12}-\lambda(1+3\mu) a_{21}a_{12}+\lambda\mu
 a_{21}a_{12}a_{21}a_{12} \\
-1-\mu+a_{21}+\mu a_{12}a_{21} &
-1-\mu+(1+2\mu)a_{12}-\mu a_{12}a_{21}a_{12}
\end{smallmatrix} \right)
\]
and
\[
A \cdot (1 - \Phir{\sigma_1^{3}}\cdot\Lambda^{-1}) = \left(
\begin{smallmatrix} 
\lambda^{-1}\left(-\lambda-\lambda\mu+(2+\mu)a_{12}
-a_{12}a_{21}a_{12}\right) &
-1-\mu+\mu a_{12}+a_{12}a_{21}\\
\lambda^{-1}\left(1+\mu+\lambda a_{21}-(3+\mu)a_{21}a_{12}+
a_{21}a_{12}a_{21}a_{12}\right) & 
-1-\mu+a_{21}+(2+\mu) a_{21}-
a_{21}a_{12}a_{21}
\end{smallmatrix} \right),
\]
once we replace $\lambda$ by $\lambda\mu^{-3}$. As in the
computation from \cite{Ng1}, it follows that in $HC_0$, we have
$a_{21}=\mu^2 a_{12}-\mu^2+1$, and, setting $a_{12}=x$, we find that
\[
HC_0(\textrm{LH trefoil}) \cong (\ring)[x]/((x-1)(\mu^2 x+\mu+1), \lambda\mu^2 x-\mu^3x-\lambda\mu^2+\lambda).
\]
We can then deduce $HC_0$ for the right handed trefoil by replacing
$\mu$ by $\mu^{-1}$:
\[
HC_0(\textrm{RH trefoil}) \cong (\ring)[x]/((x-1)(x+\mu^2+\mu), \lambda\mu x-x+\lambda\mu^3-\lambda\mu).
\]
%This agrees with the expression for the cord algebra of the right handed trefoil in Section~\ref{ssec:corddef},
%after we change variables $x\rightarrow -\mu x$.

\subsecdivide
%*********************************************************************
\subsection{Differential graded algebra interpretation II}
\label{ssec:cordDGA}

The definition of the framed knot DGA from Section~\ref{ssec:DGAdef}
is somewhat formal and inscrutable. Here we give another formulation
for the DGA, using a diagram of the knot rather than a braid
presentation. This is based on the cord formulation for the
invariant, and improves on Section~4.4 of \cite{Ng2}.

\begin{figure}
\centerline{
\includegraphics[width=3.5in]{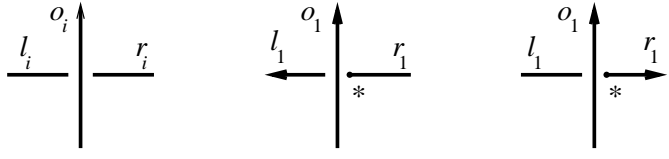}
} \caption{ Labeling the strands at a knot crossing. The left
diagram is crossing $i$ for $i\geq 2$; the center and right diagrams
are crossing $1$ for $\epsilon_1=1$ and $\epsilon_1=-1$,
respectively. } \label{fig:crossings}
\end{figure}

We first describe how to use a knot diagram to calculate the cord
algebra. Suppose that we have a knot diagram of $K$ with $n$
crossings; see Figure~\ref{fig:crossings}. The crossings divide the
diagram of the knot into $n$ connected components, which we label
$1,\ldots,n$; also number the crossings $1,\ldots,n$ in some way.
Crossing $\alpha$ can be represented as
$(o_{\alpha},l_{\alpha},r_{\alpha})$, where $o_{\alpha}$ is the
overcrossing strand, and $l_{\alpha}$ and $r_{\alpha}$ are the left
and right undercrossing strands from the point of view of the
orientation on $o_{\alpha}$. Choose the point $*$ of
Definition~\ref{def:cordcordalg} to be the undercrossing point of
crossing $1$, and define $\epsilon_1$ to be $1$ if crossing $1$ is a
positive crossing, and $-1$ if it is negative.

Let $\A_n$ denote the tensor algebra over $\ring$ generated by
$a_{ij}$, where $1\leq i,j\leq n$ and $i\neq j$, and set
$a_{ii}=1+\mu$ for $1\leq i\leq n$. Define $\I_K^{\textrm{diagram}}
\subset \A_n$ to be the subalgebra generated by the following:
\[
\begin{array}{lll}
a_{jl_i} + \mu a_{jr_i} - a_{jo_i} a_{o_il_i},~ \mu a_{l_i j} +
a_{r_i j} - a_{l_io_i} a_{o_ij},
&& 2\leq i\leq n ~\textrm{and}~ 1\leq j\leq n, \\
a_{jl_1} + \lambda^{\epsilon_1} \mu a_{jr_1} - a_{jo_1} a_{o_1l_1},~
\mu a_{l_1 j} + \lambda^{-\epsilon_1} a_{r_1 j} - a_{l_1o_1}
a_{o_1j}, && 1\leq j\leq n.
\end{array}
\]
This subalgebra is constructed to contain all of the relations
imposed by (\ref{eq:skein3new}) of Definition~\ref{def:cordcordalg};
see also \cite[\S 4.4]{Ng2}.

\begin{proposition}
The cord algebra of $(K,*)$ is given by
$\A_n/\I_K^{\textrm{diagram}}$.
\end{proposition}

\begin{proof}
Same as the proof of Proposition 4.8 in \cite{Ng2}.
\end{proof}

We can express the relations defining $\I_K^{\textrm{diagram}}$ more
neatly in terms of matrices. Write $\X_{i,j}$ for the $n\times n$
matrix which has $1$ in the $ij$ entry and $0$ everywhere else. Let
$\Psi^L$, $\Psi^R$, $\Psi^L_1$, $\Psi^R_1$, $\Psi^L_2$, $\Psi^R_2$
be the $n\times n$ matrices given as follows:
\begin{alignat*}{3}
\Psi^L_1 & = \lambda^{-\epsilon_1} \X_{1,r_1} + \sum_{\alpha\neq 1}
\X_{\alpha, r_{\alpha}} &&\hspace{0.5in}& \Psi^R_1 &=
\lambda^{\epsilon_1}\mu \X_{r_1, 1}
+ \mu \sum_{\alpha\neq 1} \X_{r_{\alpha}, \alpha} \\
\Psi^L_2 &= \sum_\alpha \left( \mu \X_{\alpha, l_{\alpha}} -
a_{l_{\alpha}o_{\alpha}} \X_{\alpha, o_{\alpha}}\right) &&& \Psi^R_2
&= \sum_\alpha \left( \X_{l_{\alpha},\alpha}
- a_{o_{\alpha}l_{\alpha}} \X_{o_{\alpha},\alpha}\right) \\
\Psi^L &= \Psi^L_1 + \Psi^L_2 &&& \Psi^R &= \Psi^R_1 + \Psi^R_2.
\end{alignat*}
This definition may seem a bit daunting, but notice that the
matrices $\Psi$ are quite sparse; for example, $\Psi^L$ (resp.\
$\Psi^R$) has at most three nonzero entries in each row (resp.\
column).

%\[
%(\Psi^L_1)_{\alpha j} =
%\begin{cases}
%1, & j=r_{\alpha} \text{ and }
%\alpha \neq 1 \\
%\lambda^{-\epsilon_1}, & j=r_{\alpha} \text{ and } \alpha=1 \\
%0, & \text{otherwise,}
%\end{cases}
%\hspace{0.5in}
%(\Psi^R_1)_{j \alpha} =
%\begin{cases}
%\mu, & j=r_{\alpha} \text{ and }
%\alpha \neq 1 \\
%\lambda^{\epsilon_1} \mu, & j=r_{\alpha} \text{ and } \alpha=1 \\
%0, & \text{otherwise,}
%\end{cases}
%\]
%\[
%(\Psi^L_2)_{\alpha j} =
%\begin{cases}
%\mu, & j=l_{\alpha}\neq o_{\alpha} \\
%-a_{l_{\alpha}o_{\alpha}}, & j=o_{\alpha}\neq l_{\alpha} \\
%\mu-a_{l_{\alpha}o_{\alpha}}, & j=l_{\alpha}=o_{\alpha} \\
%0, & \text{otherwise,}
%\end{cases}
%\hspace{0.5in}
%(\Psi^R_2)_{j\alpha} =
%\begin{cases}
%1, & j=l_{\alpha} \neq o_{\alpha} \\
%-a_{o_{\alpha}l_{\alpha}}, & j=o_{\alpha}\neq l_{\alpha} \\
%1-a_{o_{\alpha}l_{\alpha}}, & j=l_{\alpha}=o_{\alpha} \\
%0, & \text{otherwise,}
%\end{cases}
%\]
%\[
%\Psi^L = \Psi^L_1 + \Psi^L_2,
%\hspace{0.5in}
%\Psi^R = \Psi^R_1 + \Psi^R_2.
%\]

Also, assemble the generators of $\A_n$ into a matrix $\hat{A}$ with
$\hat{A}_{ii} = 1+\mu$ and $\hat{A}_{ij} = a_{ij}$ for $i\neq j$.
(Note that this is slightly different than the matrix $A$ from
Section~\ref{ssec:DGAdef}.) Then the generators of
$\I_K^{\textrm{diagram}}$ are the entries of the matrices
\[
\Psi^L \cdot \hat{A}, \hspace{0.5in} \hat{A} \cdot \Psi^R.
\]

\begin{figure}
\centerline{
\includegraphics[height=1in]{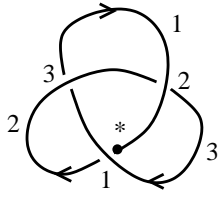}
} \caption{ A diagram of the right hand trefoil, with components
and crossings labeled. }
\label{fig:trefoil}
\end{figure}

For example, for the right handed trefoil shown in
Figure~\ref{fig:trefoil}, we have $\Psi^L = \left(
\begin{smallmatrix} \lambda^{-1} & \mu & -a_{23} \\
-a_{31} & 1 & \mu \\
\mu & -a_{12} & 1 \end{smallmatrix} \right)$ and $\Psi^R = \left(
\begin{smallmatrix} \lambda\mu & -a_{13} & 1 \\
1 & \mu & -a_{21} \\
-a_{32} & 1 & \mu
\end{smallmatrix} \right)$.
If we set the entries of $\Psi^L \cdot \hat{A}$ and $\hat{A} \cdot
\Psi^R$ equal to $0$, we find that $a_{23}=a_{21}$, $a_{32}=a_{12}$,
$a_{13}=a_{12}$, $a_{31}=a_{21}$, $\lambda a_{31} = a_{32}$,
$\lambda^{-1} a_{13} = a_{23}$. Setting $a_{12}=x$, we can use the
relations to compute that the cord ring of the right hand trefoil is
\[
(\ring)[x]/(x^2-x-\lambda\mu^2-\lambda\mu,x^2-\lambda\mu
x-\lambda-\lambda\mu).
\]
This agrees with the expression given in Section~\ref{ssec:DGAdef},
after we change variables $x\rightarrow -\mu^{-1}x$.

%As an illustration, consider the right handed trefoil $K$ shown in
%Figure~\ref{fig:trefoil}. The relations defining
%\I_K^{\textrm{diag}}$ are as follows:
%\[
%\begin{array}{ccc}
%1+\mu+\mu a_{13}-a_{12}a_{21}, & a_{21}+\mu a_{23}-(1+\mu)a_{21}, &
%a_{31}+\mu(1+\mu)-a_{32}a_{21}, \\
%\mu(1+\mu)+a_{31}-a_{12}a_{21}, & \mu a_{12}+a_{32}-(1+\mu) a_{12},
%&
%\mu a_{13}+1+\mu-a_{12}a_{23}; \\
%a_{13}+\mu a_{12}-(1+\mu) a_{13}, & a_{23}+\mu(1+\mu)-a_{21}a_{13},
%
%+\mu+\mu a_{32}-a_{31}a_{13}, \\
%mu a_{31}+a_{21}-(1+\mu)a_{31}, & \mu a_{32}+1+\mu-a_{31}a_{12}, &
%mu(1+\mu)+a_{23}-a_{31}a_{13}; \\
%a_{12}+\lambda\mu(1+\mu)-a_{13}a_{32}, & 1+\mu+\lambda\mu
%a_{21}-a_{23}a_{32}, &
%_{32}+\lambda\mu a_{31}-(1+\mu)a_{32}, \\
%mu a_{21}+(1+\mu)/\lambda-a_{23}a_{31}, &
%\mu(1+\mu)+a_{12}/\lambda-a_{23}a_{32}, & \mu
%a_{23}+a_{13}/\lambda-(1+\mu)a_{23}.
%\end{array}
%]
%Six of these relations simplify to $a_{23}=a_{21}$, $a_{32}=a_{12}$,
%$a_{13}=a_{12}$, $a_{31}=a_{21}$, $\lambda a_{31} = a_{32}$,
%$a_{13}/\lambda = a_{23}$. If we write $a_{12}=x$, then we can
%compute that the cord algebra for the right hand trefoil is
%[
%(\ring)[x]/(x^2-x-\lambda\mu^2-\lambda\mu,x^2-\lambda\mu
%x-\lambda-\lambda\mu).
%%= (\ring)[x]/((x-\mu-1)(\mu x+1),x-\lambda\mu x+\lambda\mu^2-\lambda).
%\]

We now construct a differential graded algebra whose degree $0$
homology is the cord algebra. Similarly to
Section~\ref{ssec:DGAdef}, the algebra is the graded tensor algebra
over $\ring$ generated by: $a_{ij}$, $1\leq i\neq j\leq n$, of
degree $0$; $b_{\alpha i}$ and $c_{i\alpha}$, $1\leq\alpha,i\leq n$,
of degree $1$; and $d_{\alpha\beta}$, $1\leq\alpha,\beta\leq n$, and
$e_{\alpha}$, $1\leq\alpha\leq n$, of degree $2$. The differential
is given by:
\begin{align*}
\d B &= \Psi^L \cdot \hat{A} \\
\d C &= \hat{A} \cdot \Psi^R \\
\d D &= B \cdot \Psi^R - \Psi^L \cdot C \\
\d e_i &= \left( B\cdot \Psi^R_1 - \Psi^L_2\cdot C\right)_{ii}.
\end{align*}

\begin{proposition} \label{prop:cordDGA}
This defines a differential, and the resulting differential graded
algebra is stable tame isomorphic to the framed knot DGA.
\end{proposition}

The author's current proof of Proposition~\ref{prop:cordDGA} is
quite laborious, involving verification of invariance under
Reidemeister moves and then identification with the framed knot DGA
for a particular knot diagram. Since we have no applications yet of
Proposition~\ref{prop:cordDGA}, we omit its proof. We have
nevertheless included this formulation of the framed knot DGA in
hopes that it might shed light on the topology behind the DGA
invariant.

We remark that, using the above expression for the cord algebra in
terms of a knot diagram, we can write the cord algebra for a knot of
bridge number $k$ as a quotient of $\A_k$. Express the knot as a
$2k$-plat, with $k$ vertical tangents to the knot diagram on the
left and $k$ on the right, and label the $k$ leftmost strands of the
diagram $1,\ldots,k$; then if we progressively construct the knot
starting with these $k$ strands and moving from left to right, we
can express any generator $a_{ij}$ involving a strand not among
$1,\dots,k$ in terms of the generators of $\A_k$. More generally, we
can use stable tame isomorphism to express the framed knot DGA as a
DGA all of whose generators have indices less than or equal to $k$,
although this is more difficult to establish.

\secdivide
%*********************************************************************
%*********************************************************************
\section{Equivalence and Invariance Proofs}
\label{sec:proofs}

Here we prove the main theorems from Section~\ref{sec:defs},
invariance of the framed knot DGA (Theorem~\ref{thm:DGAinv}) and
equivalence of the cord algebra and $HC_0$
(Theorem~\ref{thm:DGAcord}). The proofs rely heavily on the
corresponding proofs in \cite{Ng1,Ng2}, and the casual reader will
probably want to skip to the next section.

\subsecdivide
%*********************************************************************
\subsection{Proof of Theorem \ref{thm:DGAinv}}
\label{ssec:pfDGAinv}

Theorem~\ref{thm:DGAinv} is proven by following the proof of
Theorem~2.10 in \cite{Ng1}, which is the special case of
Theorem~\ref{thm:DGAinv} when $\lambda=\mu=1$, and making
modifications where necessary. We will provide the proof at the end
of this section.
%However, the main ideas from this proof can be seen more clearly in the proof of an intermediate
%result, Proposition~\ref{prop:HC0inv} below, which states that $HC_0$ is a framed knot
%invariant. Thus we will focus on proving Proposition~\ref{prop:HC0inv} algebraically, even though this
%is not actually necessary for the paper since it follows from the identification given by
%Theorem~\ref{thm:DGAcord} of $HC_0$ with the cord algebra.

We first need an auxiliary result, which will also be useful in
Section~\ref{ssec:symmetries}. In the definition of the framed knot
DGA, we used a matrix $\Lambda = \diag(\lambda,1,\ldots,1)$. We
could just as well have used a diagonal matrix where $\lambda$ is
the $m$-th diagonal entry rather than the first. The result states
that this does not change the framed knot DGA, up to tame
isomorphism. Note that the corresponding result for $HC_0$
essentially amounts to the statement that the cord algebra of
Definition~\ref{def:cordcordalg} is independent of the base point
$*$.

For any vector $v$ of length $n$, let $\Diag(v)$ be the $n\times n$ diagonal
matrix whose diagonal entries are the entries of $v$. If $n$ is fixed, then for
$1\leq m\leq n$, define $v_m$ to be the vector $(1,\ldots,1,\lambda,1,\ldots,1)$ of length
$n$ with $\lambda$ as its $m$-th entry. Write $\Lambda_m = \Diag(v_m)$.

\begin{proposition}
For $B\in B_n$, let $(\A,\d)$ be the framed knot DGA of $B$, as given in
Definition~\ref{def:knotDGAbraid}. Let $(\tilde{\A},\tilde{\d})$ be the DGA obtained from the
same definition, but with $\Lambda$ replaced by $\Lambda_m$ for some $m\leq n$. Then
$(\A,\d)$ and $(\tilde{\A},\tilde{\d})$ are tamely isomorphic.
\label{prop:lambda}
\end{proposition}

Before we prove Proposition~\ref{prop:lambda}, we establish a simple
lemma. Fix a braid $B\in B_n$ which closes to a knot. Let $s(B)$
denote the image of $B$ under the usual homomorphism $s:\thinspace
B_n\rightarrow S_n$ to the symmetric group; there is an obvious
action of $s(B)$ on vectors of length $n$, by permuting entries.

The matrices $\Phil{B},\Phir{B}$ have entries which are (noncommutative)
polynomials in the $a_{ij}$ variables;
for clarity, we can denote them by $\Phil{B}(A),\Phir{B}(A)$, where $A$ is the usual matrix of
the $a_{ij}$'s. If $\tilde{A}$ denotes the same matrix but with each $a_{ij}$ replaced by
$\tilde{a}_{ij}$, then we write $\Phil{B}(\tilde{A}),\Phir{B}(\tilde{A})$ to be the matrices
$\Phil{B}(A),\Phir{B}(A)$ with each $a_{ij}$ replaced by $\tilde{a}_{ij}$.

\begin{lemma}
Let $v$ be a vector of length $n$,
\label{lem:lambda}
and define $\tilde{a}_{ij}$ by the matrix equation
$\tilde{A} = \Diag(v) \cdot A \cdot \Diag(v)^{-1}$; that is,
$\tilde{a}_{ij} = \frac{v_i}{v_j} a_{ij}$. Then
\[
\Phil{B}(\tilde{A}) = \Diag(s(B)v) \cdot \Phil{B}(A) \cdot \Diag(v)^{-1}
\quad \textrm{and} \quad
\Phir{B}(\tilde{A}) = \Diag(v) \cdot \Phir{B}(A) \cdot \Diag(s(B)v)^{-1}.
\]
\end{lemma}

\begin{proof}
The $ij$ entry of $\Phil{B}(A)$ is a linear combination of words of the form
$a_{s(B)(i),i_1}a_{i_1i_2}\cdots a_{i_{k-1}i_k}a_{i_kj}$, while the $ij$ entry of
$\Phil{B}(\tilde{A})$ is the same linear combination of words
$\tilde{a}_{s(B)(i),i_1}\tilde{a}_{i_1i_2}\cdots \tilde{a}_{i_{k-1}i_k}\tilde{a}_{i_kj}
= \frac{v_{s(B)(i)}}{v_j} a_{s(B)(i),i_1}a_{i_1i_2}\cdots a_{i_{k-1}i_k}a_{i_kj}$.
Hence the $ij$ entry of $\Phil{B}(\tilde{A})$ is $\frac{v_{s(B)(i)}}{v_j}$ times the
$ij$ entry of $\Phil{B}(A)$. The result follows for $\Phil{B}$; an analogous computation
yields the result for $\Phir{B}$.
\end{proof}

\begin{proof}[Proof of Proposition~\ref{prop:lambda}]
Since the braid $B$ closes to a knot, the permutation $s(B)$ can be written as a
cycle $(1,m_1,\ldots,m_{n-1})$. Define a vector $v_m$ of length $n$ whose $i$-th entry is
$\lambda$ if $i$ appears after $1$ but not after $m$ in the cycle for $s(B)$; for example,
for $B=\sigma_1\sigma_3\sigma_2^{-3}\in B_4$, we have $s(B) = (1 2 4 3)$, and thus
$v_2 = (1,\lambda,1,1)$, $v_3=(1,\lambda,\lambda,\lambda)$, $v_4=(1,\lambda,\lambda,1)$.

To distinguish between the algebras $\A$ and $\tilde{\A}$, denote the generators of $\tilde{\A}$
by $\tilde{a}_{ij},\tilde{b}_{ij},\tilde{c}_{ij},\tilde{d}_{ij},\tilde{e}_i$; the first four
families can be assembled into matrices $\tilde{A},\tilde{B},\tilde{C},\tilde{D}$. We claim
that the identification $\tilde{A}=\Diag(v_m)\cdot A\cdot\Diag(v_m)^{-1}$,
$\tilde{B}=\Diag(v_m)\cdot B\cdot\Diag(v_m)^{-1}$,
$\tilde{C}=\Diag(v_m)\cdot C\cdot\Diag(v_m)^{-1}$,
$\tilde{D}=\Diag(v_m)\cdot D\cdot\Diag(v_m)^{-1}$,
$\tilde{e_i} = e_i$ also identifies $\tilde{\d}$ with $\tilde{\d}$.

Indeed, under this identification, we have
\begin{eqnarray*}
\tilde{\d} \tilde{B} &=& \tilde{A}-\Lambda_m \cdot \Phil{B}(\tilde{A}) \cdot \tilde{A} \\
&=& \Diag(v_m) \cdot A \cdot \Diag(v_m)^{-1}
- \Lambda_m \cdot \Diag(s(B)v_m) \cdot \Phil{B}(A) \cdot A \cdot \Diag(v_m)^{-1} \\
&=& \Diag(v_m) \cdot (A - \Lambda\cdot\Phil{B}(A)\cdot A) \cdot \Diag(v_m)^{-1} \\
&=& \d \tilde{B},
\end{eqnarray*}
where the second equality follows from \ref{lem:lambda}, and the third follows from the
fact that $\Lambda_m = \Diag(v_m) \cdot \Lambda \cdot \Diag(s(B)v_m)$, by the construction
of $v_m$. We similarly find that $\tilde{\d}=\d$ on $\tilde{C}$, $\tilde{D}$, and $\tilde{e}_i$.
Thus the identification gives a (tame) isomorphism between $\A$ and $\tilde{\A}$ which intertwines
$\d$ and $\tilde{\d}$, as desired.
\end{proof}

%If $B$ is a braid, write $HC_0^{\textrm{knot}}(B)$ for the degree $0$ homology of the framed knot DGA of $B$.
%We now proceed to give an algebraic proof of the following result, which is the analogue of
%Theorem 4.10 from \cite{Ng1}.
%
%\begin{proposition}
%If braids $B$ and $\tilde{B}$ have isotopic knot closures, then $HC_0^{\textrm{knot}}(\tilde{B})$
%is isomorphic to $HC_0^{\textrm{knot}}(B)$, after replacing $\lambda$ by $\lambda\mu^{w(\tilde{B})-w(B)}$ in the latter.
%\label{prop:HC0inv}
%\end{proposition}
%
%\begin{proof}
%It suffices to establish the proposition when $B$ and $B'$ are related by one of the Markov moves.
%We consider each in turn; our proof closely follows the proof of \cite[Thm.\ 4.10]{Ng1}, and we
%will mainly note the differences.
%
%\noindent \textsc{Conjugation.}
%Here $B,\tilde{B}$ satisfy $\tilde{B} = C^{-1}BC$ for some $C\in B_n$. The relations defining
%the degree $0$ homology of the framed knot DGA are given by the entries of
%
%\end{proof}

\begin{proof}[Sketch of proof of Theorem~\ref{thm:DGAinv}]
The theorem can be proved in a way precisely following the proof of Theorem 2.10 in \cite{Ng1}, making
alterations where necessary. It suffices to establish the theorem when $B$ and $B'$ are related by
one of the Markov moves. We consider each in turn, giving the identification of generators which
establishes the desired stable tame isomorphism, and leaving the easy but extremely tedious verifications
to the reader.

\vspace{11pt}

\noindent \textsc{Conjugation.}
Here $\tilde{B}=\sigma_k^{-1} B\sigma_k$ for some $k$. It is most convenient to prove that a slightly
different form for the framed knot DGA is the same for $B$ and $\tilde{B}$. More precisely, it is easy to check
that the framed knot DGA for any braid $B$ is stable tame isomorphic to the ``modified framed knot DGA''
with generators $\{a_{ij}\,|\,1\leq i,j\leq n, i\neq j\}$ of
degree $0$, $\{b_{ij}\,|\,1\leq i,j\leq n, i\neq j\}$ and
$\{c_{ij},d_{ij}\,|\,1\leq i,j\leq n\}$ of degree $1$,
and $\{e_{ij},f_{ij}\,|\,1\leq i,j\leq n\}$ of degree $2$, and differential
\begin{eqnarray*}
\d A &=& 0 \\
\d B &=& A - \Lambda \cdot \hom{B}(A) \cdot \Lambda^{-1} \\
\d C &=& (1 - \Lambda\cdot \Phil{B}(A)) \cdot A \\
\d D &=& A \cdot (1 - \Phir{B}(A)\cdot \Lambda^{-1}) \\
\d E &=& B - D - C \cdot \Phir{B}(A)\cdot \Lambda^{-1} \\
\d F &=& B - C - \Lambda\cdot \Phil{B}(A) \cdot D,
\end{eqnarray*}
where we set $b_{ii}=0$ for all $i$.

Define $\tilde{\Lambda}$ to be $\Lambda$ if $k\geq 2$, and
$\Lambda_2$ if $k=1$. By Proposition~\ref{prop:lambda}, the modified
framed knot DGA for $\tilde{B}$ is tamely isomorphic to the same DGA
but with $\Lambda$ replaced by $\tilde{\Lambda}$; call the latter
DGA $(\tilde{\A},\tilde{\d})$, with generators
$\tilde{a}_{ij},\tilde{b}_{ij},\tilde{c}_{ij},\tilde{d}_{ij},\tilde{e}_i$.
It suffices to exhibit a tame isomorphism between
$(\tilde{\A},\tilde{\d})$ and the modified framed knot DGA $(\A,\d)$
of $B$.

Let $\epsilon_k$ be $\lambda$ if $k=1$ and $1$ otherwise, and recall
that $a_{ij}'$ denotes the $ij$ entry in $A$, i.e., $a_{ij}'=a_{ij}$
if $i>j$ and $a_{ij}'=\mu a_{ij}$ if $i<j$. Now identify generators
of $\A$ and $\tilde{\A}$ as follows: $\tilde{a}_{ij} =
\hom{B}(a_{ij})$;
\begin{gather*}
\begin{alignat*}{4}
\tilde{b}_{ki} &= - b_{k+1,i} - \epsilon_k^{-1}
\phi(a_{k+1,k})b_{ki} - b_{k+1,k} a_{ki}' & \qquad & i\neq
k,k+1 \\
\tilde{b}_{ik} &= - b_{i,k+1} - \epsilon_k b_{ik} \phi(a_{k,k+1}) - a_{ik}' b_{k,k+1} & \qquad & i\neq
k,k+1 \\
\tilde{b}_{k+1,i} &= b_{ki} & \qquad & i\neq k,k+1 \\
\tilde{b}_{i,k+1} &= b_{ik} & \qquad & i\neq k,k+1 \\
\tilde{b}_{k,k+1} &= \mu b_{k+1,k} & \\
\tilde{b}_{k+1,k} &= \mu^{-1} b_{k,k+1} & \\
\tilde{b}_{ij} &= b_{ij} & \qquad & i,j\neq k,k+1;
\end{alignat*} \\
\tilde{C} = \Phil{\sigma_k}(\epsilon_k^{-1} \hom{B}(A)) \cdot C
\cdot \Phir{\sigma_k}(A)
+ \Theta^L_k \cdot A \cdot \Phir{\sigma_k}(A); \\
\tilde{D} = \Phil{\sigma_k}(A) \cdot D \cdot \Phir{\sigma_k}(\epsilon_k\hom{B}(A))
+ \mu^{-1} \Phil{\sigma_k}(A) \cdot A \cdot \Theta^R_k; \\
\tilde{E} = \Phil{\sigma_k}(\epsilon_k^{-1} \hom{B}(A)) \cdot E
\cdot \Phir{\sigma_k}(\epsilon_k\hom{B}(A)) - \Theta^L_k \cdot D
\cdot \Phir{\sigma_k}(\hom{B}(A))
+ \Theta^L_k \cdot \Theta^R_k; \\
\tilde{F} = \Phil{\sigma_k}(\epsilon_k^{-1} \hom{B}(A)) \cdot F
\cdot \Phir{\sigma_k}(\epsilon_k\hom{B}(A)) + \mu^{-1}
\Phil{\sigma_k}(\epsilon_k^{-1} \hom{B}(A)) \cdot C \cdot \Theta^R_k
+ \Theta^L_k \cdot (1+\mu^{-1} A) \cdot \Theta^R_k.
\end{gather*}
Here $\Phil{\sigma_k}(\epsilon_k^{-1} \hom{B}(A))$ and
$\Phir{\sigma_k}(\epsilon_k\hom{B}(A))$ are the matrices which
coincide with the $n\times n$ identity matrix except in the
intersection of rows $k,k+1$ and columns $k,k+1$, where they are
$\left( \begin{smallmatrix} -\epsilon_k^{-1}
\hom{B}(a_{k+1,k}) & -1 \\ 1 & 0 \end{smallmatrix} \right)$ and
$\left( \begin{smallmatrix} -\epsilon_k \hom{B}(a_{k,k+1}) & 1 \\ -1 & 0 \end{smallmatrix} \right)$,
respectively, and $\Theta^L_k$ and $\Theta^R_k$ are the $n\times n$ matrices which are identically
zero except in the $kk$ entry, where they are $-b_{k+1,k}$ and $-b_{k,k+1}$, respectively.

This identification of $\A$ and $\tilde{\A}$ gives the desired tame
isomorphism. For anyone wishing to verify that $\d$ and $\tilde{\d}$
are intertwined by this identification, it may be helpful to note
the following identities: $\Phil{\sigma_k}(\epsilon_k^{-1}
\hom{B}(A)) \Lambda = \tilde{\Lambda} \Phil{\sigma_k}(\hom{B}(A))$,
$\Phir{\sigma_k}(\epsilon_k\hom{B}(A)) \tilde{\Lambda} = \Lambda
\Phir{\sigma_k}(\hom{B}(A))$, and
\begin{align*}
\tilde{B} &= \Phil{\sigma_k}(\epsilon_k^{-1} \hom{B}(A)) \cdot B
\cdot \Phir{\sigma_k}(\epsilon_k \hom{B}(A)) + \mu^{-1}
\Phil{\sigma_k}(A) \cdot A
\cdot \Theta^R_k(A) \\
&\quad + \Theta^L_k(A) \cdot A \cdot
\Phir{\sigma_k}(\epsilon_k\hom{B}(A)) + \d\left( \Theta^L_k(A) \cdot
\Theta^R_k(A) \right).
\end{align*}

\vspace{11pt}

\noindent \textsc{Positive stabilization.} Here $\tilde{B}$ is
obtained by adding a strand labeled $0$ to $B$ and setting
$\tilde{B}=B \sigma_0$. We will denote the framed knot DGA of $B$ by
$(\A,\d)$ as usual. Write the generators of the framed knot DGA of
$\tilde{B}$ as
$\tilde{a}_{ij},\tilde{b}_{ij},\tilde{c}_{ij},\tilde{d}_{ij},\tilde{e}_i$,
where now $0\leq i,j\leq n$; in the definition of the differential,
the matrix $\Lambda$ now has $\lambda$ as its $(0,0)$ entry. Let
$(\tilde{\A},\tilde{\d})$ denote the result of replacing $\lambda$
by $\lambda/\mu$ in the framed knot DGA of $\tilde{B}$.

We want an identification between $\A$, suitably stabilized as in the proof of
\cite[Thm.\ 2.10]{Ng1}, and $\tilde{\A}$, which yields a stable
tame isomorphism between $(\A,\d)$ and $(\tilde{\A},\tilde{\d})$. It is given as follows:
$a_{i_1i_2} = \tilde{a}_{i_1i_2}$ for $1\leq i,j\leq n$;
\begin{xalignat*}{2}
b_{1i} &= - \mu \tilde{b}_{0i} + \tilde{b}_{1i} - \mu \tilde{b}_{00} a_{0i} & \qquad
c_{i1} &= - \mu^{-1} \tilde{c}_{i0} + \tilde{c}_{i1} - \mu^{-1} a_{i0} \tilde{c}_{00} \\
b_{ji} &= \tilde{b}_{ji} & \qquad
c_{ij} &= \tilde{c}_{ij}
\end{xalignat*}
for $i\geq 1$ and $j\geq 2$;
\begin{gather*}
\begin{split}
d_{11} &= \tilde{d}_{00} - \mu \tilde{d}_{01} + \tilde{d}_{11} -
\mu^{-1} \tilde{d}_{10} +
\tilde{b}_{00} \tilde{c}_{01} + \mu^{-1} \tilde{b}_{10} \tilde{c}_{00} \\
d_{1j} &= - \mu\tilde{d}_{0j} + \tilde{d}_{1j} + \tilde{b}_{00} \tilde{c}_{0j} \\
d_{j1} &= - \mu^{-1} \tilde{d}_{j0} + \tilde{d}_{j1} + \mu^{-1} \tilde{b}_{j0} \tilde{c}_{00} \\
d_{j_1j_2} &= \tilde{d}_{j_1j_2} \\
e_1 &= \tilde{e}_0 + \tilde{e}_1 - \mu\tilde{d}_{01} + \tilde{b}_{00} \tilde{c}_{01} \\
e_j &= \tilde{e}_j \\
d_{00} &= \mu^{-1} \tilde{d}_{10} - \tilde{d}_{11} +
\mu\tilde{d}_{01} - \mu^{-1} \tilde{b}_{10} \tilde{c}_{00}
- \tilde{b}_{00} \tilde{c}_{01} \\
e_0 &= \mu\tilde{d}_{01} - \tilde{e}_1 - \tilde{b}_{00} \tilde{c}_{01}
\end{split} \\
d_{10} = \tilde{d}_{11} - \tilde{e}_1 \qquad
d_{01} = \tilde{e}_1 \qquad
d_{0j} = \tilde{d}_{1j} \qquad
d_{j0} = \tilde{d}_{j1}
\end{gather*}
for $j,j_1,j_2\geq 2$; and
\begin{xalignat*}{3}
b_{10} &= \tilde{b}_{10} + \tilde{c}_{11} & \textrm{and} &&
c_{01} &= \tilde{c}_{01} + \tilde{b}_{11}; \\
b_{00} &= \tilde{b}_{00} - b_{11} - \lambda\Phi^L_{1\ell} \tilde{c}_{\ell 1} & \textrm{and} &&
c_{00} &= \tilde{c}_{00} - c_{11} - \lambda^{-1} \tilde{b}_{1\ell} \Phi^R_{\ell 1}; \\
b_{j0} &= \tilde{b}_{j0} + \tilde{c}_{j1} - b_{j1} - \Phi^L_{j\ell} \tilde{c}_{\ell 1} & \textrm{and} &&
c_{0j} &= \tilde{c}_{0j} + \tilde{b}_{1j} - c_{1j} - \tilde{b}_{1\ell}\Phi^R_{\ell j}; \\
b_{0i} &= \tilde{b}_{1i} & \textrm{and} &&
c_{i0} &= \tilde{c}_{i1}
\end{xalignat*}
for $i\geq 1$ and $j \geq 2$. If we set $\d b_{0i} = -\mu a_{0i}+a_{1i}'$,
$\d c_{i0} = -a_{i0}+a_{i1}'$, $\d e_0 = b_{00}$, $\d d_{00} = b_{00} - c_{00}$, $\d d_{i0} = -b_{i0}$,
$\d d_{0i} = c_{0i}$ for $i\geq 1$, then the above identification produces $\d=\tilde{\d}$.

\vspace{11pt}

\noindent \textsc{Negative stabilization.}
Here $\tilde{B}=B \sigma_0^{-1}$; let $(\A,\d)$ be the framed knot DGA of $B$, and let $(\tilde{\A},\tilde{\d})$
be the framed knot DGA of $\tilde{B}$ with $\lambda$ replaced by $\lambda\mu$, with notation as in
the previous case.

In this case, the identification between $\A$ and $\tilde{\A}$ is as follows:
$a_{i_1i_2} = \tilde{a}_{i_1i_2}$ for $1\leq i,j\leq n$;
\begin{xalignat*}{2}
b_{1i} &= \tilde{b}_{1i} - \tilde{b}_{0i} + \lambda \tilde{c}_{00}
\Phi^L_{1\ell} a_{\ell i} & \qquad
c_{i1} &= \tilde{c}_{i1} - \tilde{c}_{i0} + \textstyle{\frac{1}{\lambda\mu}} a_{i\ell} \Phi^R_{\ell 1} \tilde{b}_{00} \\
b_{ji} &= \tilde{b}_{ji} & \qquad
c_{ij} &= \tilde{c}_{ij};
\end{xalignat*}
\begin{gather*}
\begin{split}
d_{11} &= - \tilde{d}_{10} + \tilde{d}_{11} - \tilde{d}_{01} + \tilde{d}_{00}
- \textstyle{\frac{1}{\lambda\mu}} b_{1\ell} \Phi^R_{\ell 1} \tilde{b}_{00} - \lambda \tilde{c}_{00} \Phi^L_{1\ell} c_{\ell 1}
-\textstyle{\frac{\mu+1}{\mu}} \tilde{c}_{00}\tilde{b}_{00} \\
d_{1j} &= \tilde{d}_{1j} - \tilde{d}_{0j} - \lambda \tilde{c}_{00} \Phi^L_{1\ell} c_{\ell j} \\
d_{j1} &= \tilde{d}_{j1} - \tilde{d}_{j0} - \textstyle{\frac{1}{\lambda\mu}} b_{j\ell} \Phi^R_{\ell 1} \tilde{b}_{00} \\
d_{j_1j_2} &= \tilde{d}_{j_1j_2} \\
e_1 &= -\tilde{d}_{01} + \tilde{e}_0 + \tilde{e}_1 - \lambda \tilde{c}_{00} \Phi^L_{1\ell} c_{\ell 1}
- \tilde{c}_{00} \tilde{b}_{00} + \lambda\mu (\tilde{e}_0-\tilde{d}_{00})\hom{B}(a_{10})
+ \textstyle{\frac{1}{\lambda\mu}} \hom{B}(a_{01})\tilde{e}_0 \\
e_j &= \tilde{e}_j \\
\end{split} \\
d_{00} = \tilde{d}_{00} \qquad
e_0 = \tilde{e}_0 \qquad
d_{10} = \tilde{d}_{11} - \tilde{e}_1 \qquad
d_{01} = \tilde{e}_1 \qquad
d_{j0} = \tilde{d}_{j0} \qquad
d_{0j} = \tilde{d}_{0j};
\end{gather*}
\begin{xalignat*}{3}
\begin{split}
b_{10}
&= -\tilde{b}_{10} + c_{11} + \tilde{c}_{10} - \textstyle{\frac{1}{\lambda\mu}} a_{1\ell}' \Phi^R_{\ell 1} \tilde{b}_{00} \\
& \quad - \{b_{1\ell}\Phi^R_{\ell 1} +
\tilde{b}_{0\ell}\Phi^R_{\ell 1}\\
& \quad\quad + \lambda(1+\mu) \tilde{c}_{00}\} \hom{B}(a_{10})
\end{split}
& \textrm{and} &&
\begin{split}
c_{01} &= -\tilde{c}_{01} + b_{11} + \tilde{b}_{01} - \lambda \tilde{c}_{00}\Phi^L_{1\ell} a_{\ell 1}' \\
& \quad - \hom{B}(a_{01}) \{\Phi^L_{1\ell} c_{\ell 1} \\
& \quad\quad + \Phi^L_{1\ell} \tilde{c}_{\ell 0}  +
\textstyle{\frac{1+\mu}{\lambda\mu}} \tilde{b}_{00}\};
\end{split} \\
b_{00} &= \tilde{b}_{00} + \lambda\mu \Phi^L_{1\ell} \tilde{c}_{\ell 0} & \textrm{and} &&
c_{00} &= \tilde{c}_{00} + \textstyle{\frac{1}{\lambda\mu}} \tilde{b}_{0\ell} \Phi^R_{\ell 1}; \\
b_{j0} &= - \tilde{b}_{j0} + \tilde{c}_{j0}
+ \textstyle{\frac{1}{\lambda\mu}} b_{j\ell} \Phi^R_{\ell 1} - \Phi^L_{j\ell} \tilde{c}_{\ell 0} & \textrm{and} &&
c_{0j} &= - \tilde{c}_{0j} + \tilde{b}_{0j}
+ \lambda\mu \Phi^L_{1\ell} c_{\ell j} - \tilde{b}_{0\ell} \Phi^R_{\ell j}; \\
b_{0i} &= \tilde{b}_{0i} & \textrm{and} &&
c_{i0} &= \tilde{c}_{i0}.
\end{xalignat*}
Then $\d=\tilde{\d}$ as before, once we set
$\d b_{0i} = \mu a_{0i} - \lambda\mu \Phi^L_{1\ell} a_{\ell i}$,
$\d c_{i0} = a_{i0} - \frac{1}{\lambda\mu} a_{i\ell} \Phi^R_{\ell 1}$,
$\d d_{00} = b_{00}-c_{00}$, $\d e_0 = b_{00}$,
$\d d_{i0} = - b_{i0}$, $\d d_{0i} = c_{0i}$.
\end{proof}

\subsecdivide
%*********************************************************************
\subsection{Proof of Theorem \ref{thm:DGAcord}}
\label{ssec:pfDGAcord}

This proof is essentially the proof of Theorem 1.3 from \cite{Ng2},
with only slight modifications. Rather than giving the proof in
full, we will outline the steps here and refer the reader to
\cite{Ng2} for details.

\begin{figure}
\centerline{
\includegraphics[width=4in]{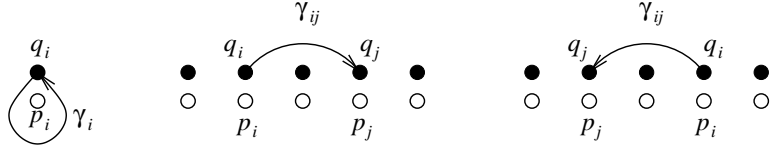}
}
\caption{
The framed arcs $\gamma_i$ (left) and $\gamma_{ij}$ for
$i<j$ (middle) and $i>j$ (right).
}
\label{fig:gammaij}
\end{figure}

Let $D$ be the unit disk in $\C$, and choose points
$p_1,\dots,p_n,q_1,\dots,q_n$ in $D$ such that $p_1,\dots,p_n$ lie
on the real line with $p_1<\dots<p_n$ and $q_j = p_j+\epsilon i$ for
each $j$ and some small $\epsilon>0$. We will depict the points
$p_1,\dots,p_n$ by an open circle, and $q_1,\dots,q_n$ by closed
circles; see Figure~\ref{fig:gammaij}.

The braid group $B_n$ acts by homeomorphisms in the usual way on the
punctured disk $D\setminus P$, where $P$ denotes
$\{p_1,\dots,p_n\}$. That is, $\sigma_k$ rotates $p_k,p_{k+1}$ in a
counterclockwise fashion around each other until they change places.
Furthermore, we can stipulate that the line segments $p_kq_k$ and
$p_{k+1}q_{k+1}$ remain vertical and of fixed length throughout the
isotopy, so that $\sigma_k$ also interchanges $q_k,q_{k+1}$, and
hence every element of $B_n$ fixes the set $Q=\{q_1,\dots,q_n\}$.
Then $B_n$ acts on the set of framed cords, as defined below.

\begin{definition}[{cf.\ \cite[Def. 3.1]{Ng2}}] \label{def:diskcord}
A \textit{framed cord} of $(D,P,Q)$ is a continuous map
$\gamma:\thinspace [0,1]\rightarrow D\setminus P$ such that
$\gamma(0),\gamma(1)\in Q$. The set of framed arcs modulo homotopy
is denoted $\P_n$. Define $\delta_i$, $1\leq i\leq n$, and
$\gamma_{ij}$, $1\leq i\neq j\leq n$, to be the framed cords
depicted in Figure~\ref{fig:gammaij}.
\end{definition}

\begin{proposition}[{cf.\ \cite[Props.\ 2.2, 3.2, 3.3]{Ng2}}] \label{prop:cordmap}
With $\A_n$ as in Section~\ref{ssec:DGAdef}, there is an unique map
$\psi:\thinspace  \P_n\rightarrow \A_n$ satisfying:
\begin{enumerate}
\item \label{item:equivar} Equivariance: $\psi(B\cdot\gamma) = \phi_B(\psi(\gamma))$ for
any $B\in B_n$ and $\gamma\in\P_n$, where $B\cdot\gamma$ denotes the
action of $B$ on $\gamma$;
\item \label{item:fram} Framing: if $\gamma$ is a framed arc with $\gamma(0)=p_j$,
then $\psi(\gamma_j\gamma) = \mu\psi(\gamma)$, where
$\gamma_j\gamma$ denotes the framed arc obtained by concatenating
$\gamma_j$ and $\gamma$m and similarly if $\gamma(1)=p_l$, then
$\psi(\gamma\gamma_l) = \mu\psi(\gamma)$;
\item \label{item:norm} Normalization: $\psi(\gamma_{ij})=-a_{ij}$ if $i<j$,
$\psi(\gamma_{ij})=-\mu a_{ij}$ if $i>j$, and $\psi(\gamma_i^0) = 1+\mu$ where
$\gamma_i^0$ is the trivial cord which is constant at $q_i$.
\end{enumerate}
Furthermore, we have the following skein relations for framed cords:
\begin{gather*}
\psi(\raisebox{-0.12in}{\includegraphics[height=0.35in]{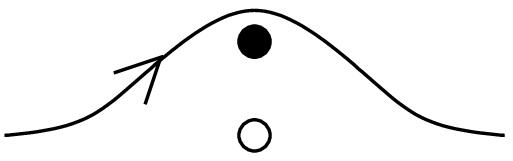}})
+
\psi(\raisebox{-0.12in}{\includegraphics[height=0.35in]{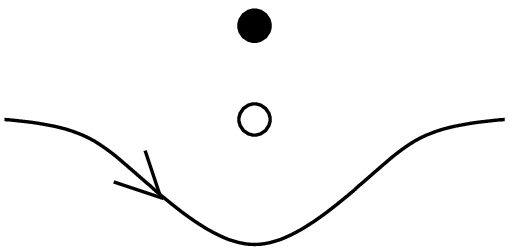}})
=
\psi(\raisebox{-0.12in}{\includegraphics[height=0.35in]{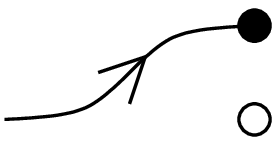}})
\psi(\raisebox{-0.12in}{\includegraphics[height=0.35in]{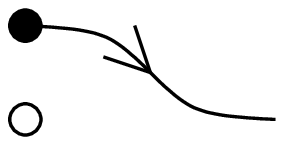}}),\\
\psi(\raisebox{-0.12in}{\includegraphics[height=0.35in]{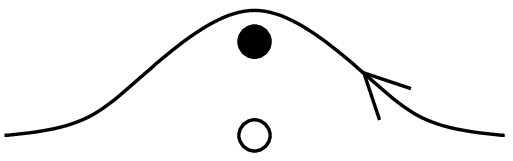}})
+
\psi(\raisebox{-0.12in}{\includegraphics[height=0.35in]{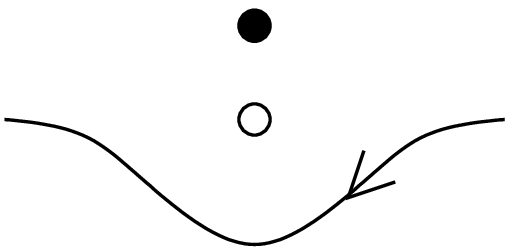}})
=
\psi(\raisebox{-0.12in}{\includegraphics[height=0.35in]{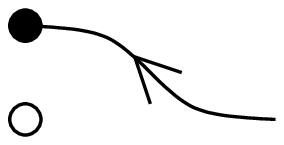}})
\psi(\raisebox{-0.12in}{\includegraphics[height=0.35in]{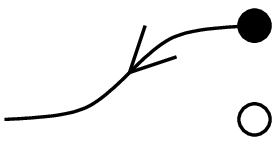}}).
\end{gather*}
Finally, the above skein relations, along with (\ref{item:fram}) and
(\ref{item:norm}), suffice to define $\psi$.
\end{proposition}

\noindent Note that the normalization of (\ref{item:norm}) has been
chosen for compatibility with (\ref{item:equivar}) and
(\ref{item:fram}): $\sigma_k$ sends $\gamma_{k,k+1}$ to
\raisebox{-0.15in}{\includegraphics[height=0.3in]{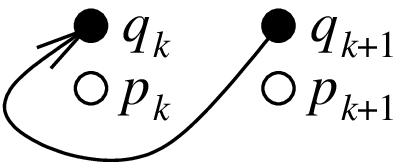}},
which is mapped under $\psi$ to $\mu^{-1}\psi(\gamma_{k+1,k})$. The
main difference between Proposition~\ref{prop:cordmap} and the
corresponding results from \cite{Ng2} is the framing axiom, which
incorporates the variable $\mu$ into the map $\psi$.

\begin{proof}[Proof of Theorem~\ref{thm:DGAcord} (cf.\ {\cite[\S 3.2]{Ng2}})]
Suppose that the knot $K$ is the closure of a braid $B\in B_n$, with framing
induced from the braid. If $\ell$ is a line
in $\R^3$, then $\R^3\setminus\ell$ is topologically a solid torus $S^1\times D^2$;
let $D$ be a $D^2$ fiber of this solid torus, viewed as sitting in $\R^3$. The
complement of $D$ in $S^1\times D^2$ is a solid cylinder in which we can embed
$B$, and the closure of this embedding in $\R^3$ is precisely $K$.

\begin{figure}
\centerline{
\includegraphics[height=2in]{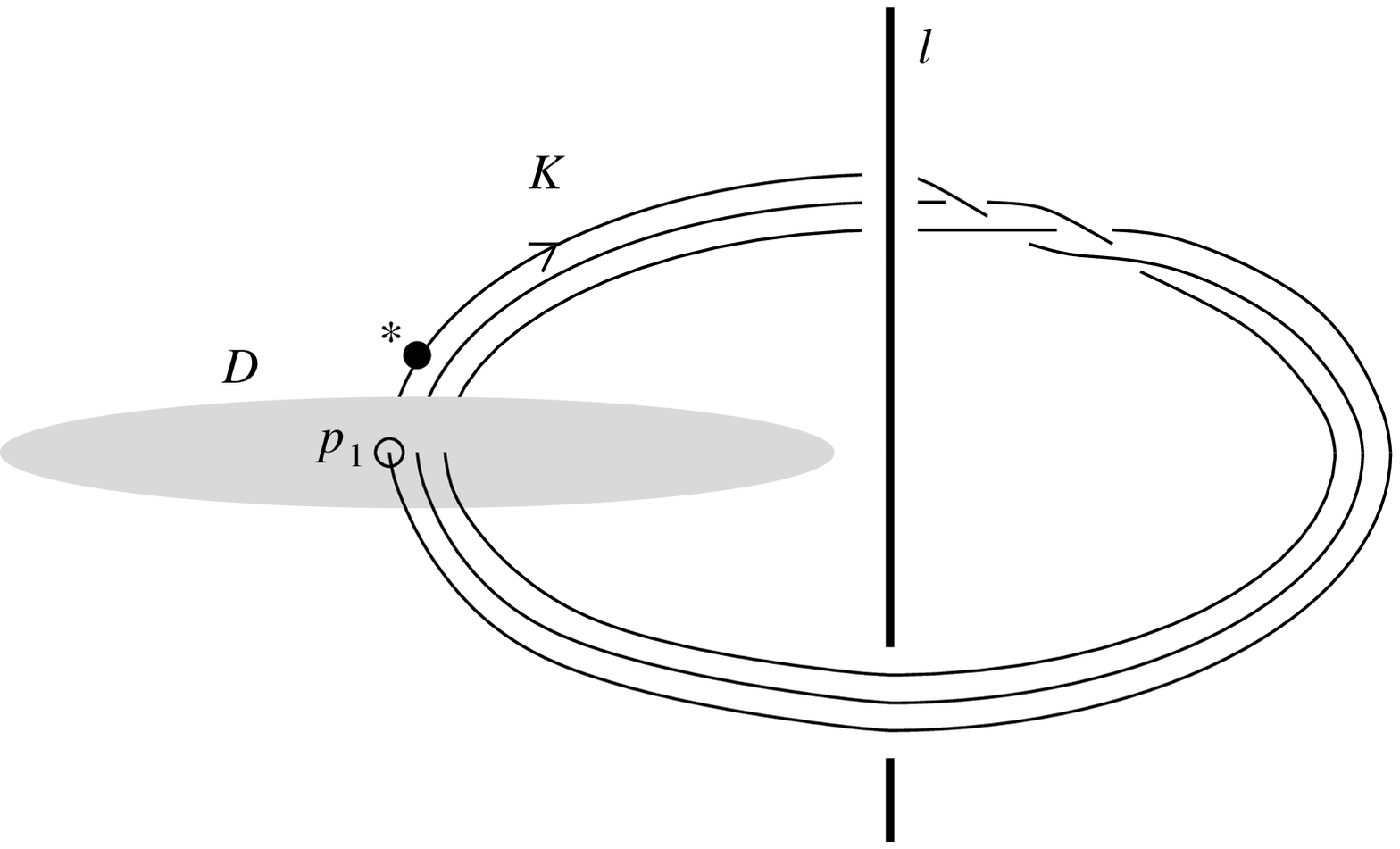}
\hspace{0.5in}
\includegraphics[height=2in]{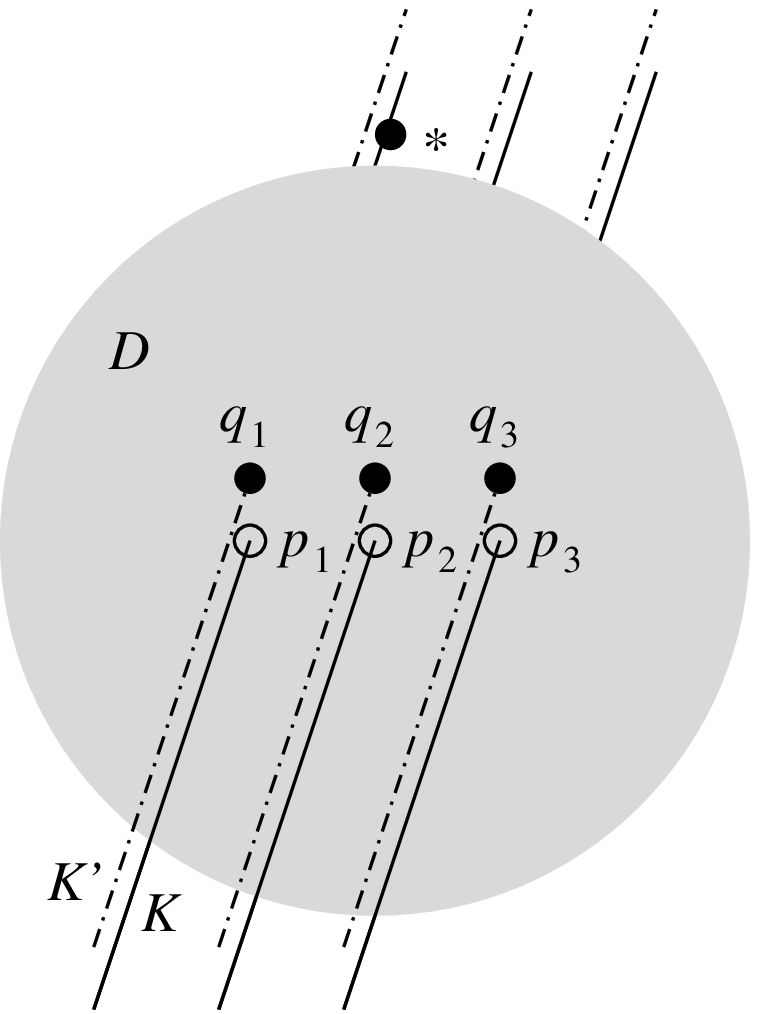}
}
\caption{Wrapping a braid around the line $\ell$ to create the
knot $K$, and a close-up of the area around the disk $D$.}
\label{fig:knotbraidpf}
\end{figure}

Now $K$ intersects $D$ in $n$ points $p_1,\dots,p_n$, and the
framing of $K$ allows us to push $K$ off of itself to obtain another
knot $K'$, which then intersects $D$ in points $q_1,\dots,q_n$.
Choose the base point $*$ on $K$, which we need to define the cord
algebra, to be just beyond $p_1$ in the direction of the orientation
of $K$. See Figure~\ref{fig:knotbraidpf}. For $1\leq i\leq n$, let
$\eta_i$ denote the path along $K'$ running from $q_i$ to $q_1$ in
the direction of the orientation of $K$ (so that it does not pass
$*$).

Any cord $\gamma$ of $K$, in the sense of Definition~\ref{def:cord},
can be homotoped to lie entirely in $D$; a slight perturbation of
the ends yields a framed cord $\tilde{\gamma}$ of $(D,P,Q)$, in the
sense of Definition~\ref{def:diskcord}, say beginning at $q_i$ and
ending at $q_j$. Now the perturbation yielding $\tilde{\gamma}$ is
not uniquely defined up to homotopy, and so $\psi(\tilde{\gamma})$
is defined in $\A_n$ only up to multiplication by powers of $\mu$.
We can however compensate for this by considering the linking number
of $K$ and the closed oriented curve $\tilde{\gamma} \cup \eta_j
\cup (-\eta_i)$. More precisely, it is easy to check that
$\mu^{\lk(K,\tilde{\gamma} \cup \eta_j \cup (-\eta_i))}
\psi(\tilde{\gamma})$ is well-defined. Given a cord of $K$ lying in
$D$, we hence get an element of $\A_n$ which we denote by
$\tilde{\psi}(\gamma)$. Moreover, the skein relation defining the
cord algebra maps to the skein relation in
Proposition~\ref{prop:cordmap}.

We conclude that the cord algebra of $K$ is a quotient of $\A_n$;
this quotient arises from the non-uniqueness of the homotopy from an
arbitrary cord of $K$ to a cord of $K$ lying in $D$. Now the
argument of the proof of \cite[Thm.\ 1.3]{Ng2} shows that the cord
algebra of $K$ is $\A_n/\I$, where $\I$ is the subalgebra of $\A_n$
generated by the entries of the matrices
$(1-\Lambda\cdot\Phil{B})\cdot A$ and $A\cdot
(1-\Phir{B}\cdot\Lambda^{-1})$. (The reason that we need the entries
of $A$ above the diagonal to be of the form $\mu a_{ij}$, rather
than $a_{ij}$ as in \cite{Ng2}, is because in this case the map
$\psi$ sends $\gamma_{ij}$ to $-\mu a_{ij}$ rather than $-a_{ij}$.)
But $\A_n/\I$ is precisely the degree $0$ homology of the framed
knot DGA of $K$.
\end{proof}

\secdivide
%*********************************************************************
%*********************************************************************
\section{Properties of the Invariants}
\label{sec:properties}

%*********************************************************************
\subsection{Symmetries}
\label{ssec:symmetries}
The cord algebra and framed knot DGA
invariants change in predictable manners when we change the framed
knot by various elementary operations: framing change, mirror, and
inversion. Here, as usual, ``inversion'' refers to reversing the
orientation of the knot, and the framing under mirror and
inversion is the induced framing from the original knot.

The effects of these operations are easiest to see in the homotopy
interpretation of the framed cord algebra
(Definition~\ref{def:knotcordalg}). As mentioned previously,
changing the framing of a knot by $f$ keeps the meridian fixed and
changes the longitude by the corresponding power of the meridian;
hence the framed cord algebra changes by replacing $\lambda$ by
$\lambda\mu^f$. A reflection of three-space about a plane changes a
framed knot into its mirror, maps the knot's longitude to its
mirror's longitude, and maps the meridian to the mirror's meridian
with the opposite orientation; hence mirroring changes the framed
cord algebra by replacing $\mu$ by $\mu^{-1}$. Changing the
orientation of a framed knot changes the orientation of both its
longitude and meridian, and so inversion replaces $\lambda$ by
$\lambda^{-1}$ and $\mu$ by $\mu^{-1}$.

\begin{proposition}
The framed cord algebra of a framed knot changes under the following
operations as shown: framing change by $f$, $\lambda\rightarrow
\lambda\mu^f$; mirror, $\mu\rightarrow\mu^{-1}$; inversion,
$\lambda\rightarrow\lambda^{-1}$ and $\mu\rightarrow\mu^{-1}$. In
addition, the commutative framed cord algebras (in which
multiplication is commutative) of a knot and its inverse are
isomorphic.
\label{prop:cordalgproperties}
\end{proposition}

\begin{proof}
The only statement left to prove is the invariance of the
commutative framed cord algebra under inversion. Let $l,m$ be the
longitude and meridian of a framed knot $K$; then $l^{-1},m^{-1}$
are the longitude and meridian of its inverse. The map sending
$\gamma\in\pi_1(S^3\setminus K)$ to $\gamma^{-1}$ extends in a
natural way to a map on the tensor algebra generated by
$\pi_1(S^3\setminus K)$. This map descends to an isomorphism on the
commutative framed cord algebras; the key point is that the relation
$[\gamma_1\gamma_2] + [\gamma_1m\gamma_2] = [\gamma_1\gamma_2]$ in
the cord algebra of $K$ is mapped to
\[
[\gamma_2^{-1}\gamma_1^{-1}] + [\gamma_2^{-1}m^{-1}\gamma_1^{-1}]
= [\gamma_1]^{-1}[\gamma_2]^{-1} = [\gamma_2]^{-1}[\gamma_1]^{-1},
\]
which is the corresponding relation for the cord algebra of the inverse of $K$.
\end{proof}

Not surprisingly, the same operations lead to the same
transformations in $\lambda,\mu$ in the framed knot DGA, although
the proof is not particularly enlightening. Define the
\textit{commutative framed knot DGA} of a knot to be the quotient of
the framed knot DGA obtained by imposing sign-commutativity, i.e.,
$vw=(-1)^{|v||w|}wv$. (Note that this does \textit{not} specialize
to the ``abelian knot DGA'' from \cite{Ng1} when we set
$\lambda=\mu=1$; there does not seem to be a reasonable lift of the
abelian knot DGA to group ring coefficients.)

\begin{proposition}
The framed knot DGA of a framed knot changes under framing change, mirroring, and inversion
as in Proposition~\ref{prop:cordalgproperties}, and the commutative framed knot DGAs of a knot
and its inverse are equivalent.
\label{prop:DGAproperties}
\end{proposition}

\begin{proof}
Framing change has already been covered by Theorem~\ref{thm:DGAinv}.

\textit{Mirror:} We follow the proof of Proposition 6.9 from
\cite{Ng1}, which is the analogous result for the DGA over $\Z$; see
that proof for notation. If a framed knot $K$ is given as the
closure of a braid $B\in B_n$, then the mirror of $K$ (with the
corresponding framing) is the closure of $B^*$. Let $(\A,\d)$ be the
framed knot DGA for $B$, and let $(\tilde{\A},\tilde{\d})$ be the
framed knot DGA for $B^*$, but with $\mu$ replaced by $\mu^{-1}$ and
$\Lambda$ replaced by $\tilde{\Lambda} = \Lambda_n$, with notation
as in Section~\ref{ssec:pfDGAinv}; this is equivalent to the usual
framed knot DGA for $B^*$ (with $\mu\rightarrow\mu^{-1}$) by
Proposition~\ref{prop:lambda}. It suffices to show that $(\A,\d)$
and $(\tilde{\A},\tilde{\d})$ are tamely isomorphic.

Assemble the generators $a_{ij},\tilde{a}_{ij}$ into matrices
$A,\tilde{A}$ as usual, except that $\tilde{A}$ has $\mu$ replaced
by $\mu^{-1}$. If we set $\tilde{a}_{n+1-i,n+1-j}=a_{ij}$ as in the
proof of \cite[Prop.\ 6.9]{Ng1}, then we find that
$\tilde{A}=\mu^{-1}\Xi(A)$, while $\Phil{B^*}(\tilde{A}) =
\Xi(\Phil{B}(A))$ and $\Phir{B^*}(\tilde{A}) = \Xi(\Phir{B}(A))$.
Now identify the other generators of $\A$ and $\tilde{\A}$ by
setting $\tilde{B}=\mu^{-1}\Xi(B)$, $\tilde{C}=\mu^{-1}\Xi(C)$,
$\tilde{D}=\mu^{-1}\Xi(D)$, $\tilde{e}_i=\mu^{-1} e_{n+1-i}$. Then
\[
\d\tilde{B} = \mu^{-1}\d\Xi B = \mu^{-1}\Xi((1-\Lambda\cdot\Phil{B}(A))\cdot A)
= (1-\tilde{\Lambda}\cdot \Phil{B^*}(\tilde{A}))\cdot\tilde{A}
= \tilde{\d}\tilde{B}
\]
and similarly $\d=\tilde{\d}$ on the other generators of $\A$.

\textit{Inverse:} Let $(\A,\d)$ be the framed knot DGA for $B$, and
let $(\tilde{\A},\tilde{\d})$ be the framed knot DGA for $B^{-1}$,
but with $\lambda$ replaced by $\lambda^{-1}$. Since the closure of
$(B^{-1})^*$ is the inverse of the closure of $B$, and the framed
knot DGA for $(B^{-1})^*$ is the framed knot DGA for $B^{-1}$ with
$\mu$ replaced by $\mu^{-1}$ by the above proof for mirrors, it
suffices to prove that $(\A,\d)$ and $(\tilde{\A},\tilde{\d})$ are
equivalent to conclude that the framed knot DGAs of a knot is
equivalent to the framed knot DGA of its inverse after setting
$\lambda\rightarrow\lambda^{-1}$ and $\mu\rightarrow\mu^{-1}$.

Assemble the generators of $\tilde{\A}$ into matrices $\tilde{A}$, etc.,
and let $\tilde{\Lambda}=\Lambda^{-1}$, so that
$\tilde{\d}\tilde{B} = (1-\tilde{\Lambda}\cdot\Phil{B^{-1}}(\tilde{A}))\cdot\tilde{A}$
and so forth. Now set $\tilde{a}_{ij}=\hom{B}(a_{ij})$, which implies that
$\tilde{A} = \Phil{B}(A)\cdot A\cdot\Phir{B}(A)$, and
$\tilde{B} = -\Lambda^{-1}\cdot B\cdot\Phir{B}(A)$,
$\tilde{C} = -\Phil{B}(A)\cdot C\cdot\Lambda$,
$\tilde{D} = \Lambda^{-1}\cdot D\cdot\Lambda$,
$\tilde{e}_i = d_{ii}-e_i$.
Since $\Phil{B^{-1}}(\tilde{A}) =(\Phil{B}(A))^{-1}$ and
$\Phir{B^{-1}}(\tilde{A})=(\Phir{B}(A))^{-1}$, it is easy to check that
this identification gives $\d=\tilde{\d}$; for instance,
\[
\tilde{\d}\tilde{B} = (1-\Lambda^{-1}\cdot(\Phil{B}(A))^{-1})\cdot\Phil{B}(A)\cdot A\cdot\Phir{B}(A)
= -\Lambda^{-1}\cdot(1-\Lambda\cdot\Phil{B})\cdot A\cdot\Phir{B}(A) = \d\tilde{B}.
\]

It remains to show that the commutative framed knot DGAs of a knot
and its inverse are equivalent. Because of the previous argument, it
suffices to prove that the commutative framed knot DGA $(\A,\d)$ of
a braid $B$ is equivalent to the same commutative framed knot DGA
$(\tilde{\A},\tilde{\d})$ with $\lambda,\mu$ replaced by
$\lambda^{-1},\mu^{-1}$. As usual, assemble the generators of
$\tilde{\A}$ into matrices $\tilde{A}$, etc., except that
$\tilde{A}$ has $\mu$ replaced by $\mu^{-1}$. If we set
$\tilde{a}_{ij} = a_{ji}$ for all $i,j$, then $\tilde{A} = \mu^{-1}
A^t$. Further set $\tilde{B}=\mu^{-1}C^t$, $\tilde{C}=\mu^{-1}B^t$,
$\tilde{D}=-\mu^{-1}D^t$, $e_i=\mu^{-1}(e_i-d_{ii})$. Now since the
$a_{ij}$ commute with each other, Proposition 6.1 from \cite{Ng1}
implies that $\Phil{B}(\tilde{A}) = (\Phir{B}(A))^t$ and
$\Phir{B}(\tilde{A}) = (\Phil{B}(A))^t$. It follows easily that
$\d=\tilde{\d}$; for instance,
\[
\tilde{\d}\tilde{B} = (1-\Lambda^{-1}\cdot\Phil{B}(\tilde{A}))\cdot\tilde{A}
= \mu^{-1}(1-\Lambda^{-1}\cdot(\Phir{B}(A))^t)\cdot A^t
= \mu^{-1}\left(A\cdot(1-\Phir{B}(A))\right)^t = \d \tilde{B},
\]
as claimed.
\end{proof}

Although the commutative framed knot DGA does not distinguish
between knots and their inverses, it is conceivable that the usual
noncommutative framed knot DGA, or the cord algebra, might.
Presumably an application to distinguish inverses would require some
sort of involved calculation using noncommutative Gr\"obner bases.

\subsecdivide
%*********************************************************************
\subsection{Relation to the alternate knot DGA}
\label{ssec:ZDGA}

It is immediate from the definition of the framed knot DGA that one can recover
the knot DGA over $\Z$ from \cite{Ng1} by setting $\lambda=\mu=1$.
What is less obvious is that the framed knot DGA also specializes to the
``alternate knot DGA'' from \cite[\S 9]{Ng1}.

\begin{proposition}
The knot DGA is obtained from the framed knot DGA by setting
$\lambda=\mu=1$; the alternate knot DGA is obtained from the framed
knot DGA by setting $\lambda=\mu=-1$.
%$\mu=-1$, $\lambda=(-1)^{f+1}$, where $f$ is the framing of the knot.
\label{prop:altDGA}
\end{proposition}

\begin{proof}
Let $B\in B_n$ be a braid whose closure gives the desired framed
knot, and note that, since $B$ closes to a knot, $w(B)\equiv n+1
\pmod 2$. Setting $\lambda=\mu=-1$ in the framed knot DGA for $K$ is
the same as setting $\mu=-1$, $\lambda=(-1)^{w(B)+1}=(-1)^n$ in the
framed knot DGA corresponding to $B$. The argument of the proof of
Proposition~\ref{prop:lambda} shows that $\Lambda$ can be replaced
in the definition of the framed knot DGA by any diagonal matrix, the
product of whose diagonal entries is $\lambda$. Since
$\lambda=(-1)^n$, we can use the matrix $-1$ instead of $\Lambda$ in
defining the framed knot DGA. Let $(\A,\d)$ denote the resulting
DGA, let $(\tilde{\A},\tilde{\d})$ denote the alternate knot DGA of
$B$, and assemble the generators of $\A$ and $\tilde{\A}$ into
matrices, as usual.

If we set $a_{ij}=\tilde{a}_{ij}$ if $i<j$ and $a_{ij}=-\tilde{a}_{ij}$
if $i>j$, then $A=-\tilde{A}$. In the notation of Section 9.1 of \cite{Ng1},
since $\phi$, $\tilde{\phi}$ are conjugate through the automorphism which
negates $a_{ij}$ for $i>j$, we have $\tilde{\Phi}^L_B(\tilde{A}) = \Phil{B}(A)$
and $\tilde{\Phi}^R_B(\tilde{A})=\Phir{B}(A)$. It follows readily from the
definition of the alternate knot DGA that if
we further set $\tilde{B}=-B$, $\tilde{C}=-C$, $\tilde{D}=-D$, $\tilde{e}_i=-e_i$,
then $\d=\tilde{\d}$. The proposition follows.
\end{proof}

Note that there are, in fact, four natural projections of the framed knot
DGA to an invariant over $\Z$, corresponding to sending each of $\lambda,\mu$ to
$\pm 1$. We have seen that two of these projections yield the knot DGA and the
alternate knot DGA from \cite{Ng1}. The other two seem to behave slightly less nicely.
For instance, for two-bridge knots, it can be shown that both the knot DGA and
the alternate knot DGA have degree $0$ homology which can be written as a quotient
of $\Z[x]$ by a principal ideal. (For the knot DGA, this was shown in
\cite[Thm.\ 7.1]{Ng1}.) This is not the case for the other two projections.

As an example, consider the left-handed trefoil examined in
Section~\ref{ssec:DGAdef}. For $\lambda=\mu=1$, we obtain $HC_0
\cong \Z[x]/(x^2+x-2)$, the knot DGA for the trefoil, and for
$\lambda=\mu=-1$, $HC_0 \cong \Z[x]/(x^2-x)$, the alternate knot DGA
for the trefoil. By contrast, we have $HC_0 \cong \Z[x]/(x^2-x,2x)$
when $\lambda=1$ and $\mu=-1$, and $HC_0 \cong \Z[x]/(x^2+x-2,2x)$
when $\lambda=-1$ and $\mu=1$.

It may be worth noting that, by
Proposition~\ref{prop:DGAproperties}, none of the projections to
$\Z$ can detect mirrors or inverses, since $\lambda=\lambda^{-1}$
and $\mu=\mu^{-1}$ in all cases. We also remark that the ``special''
status of the cases $\lambda=\mu=1$ and $\lambda=\mu=-1$ is an
artifact of a choice of a particular spin structure on $T^2$. More
precisely, the signs in the definition of the framed knot DGA, which
arise from a coherent choice of orientations of the relevant moduli
spaces, depend on a choice of one of the (four) spin structures on
the Legendrian $2$-torus; changing the spin structure negates either
or both of $\lambda$ and $\mu$. See \cite{EES2}.

\subsecdivide
%*********************************************************************
\subsection{Framed knot contact homology and the Alexander
polynomial}

It has been demonstrated in \cite[\S 7.2]{Ng2} that a linearized
version of knot contact homology over $\Z$ encodes the determinant
of the knot. Here we show that the framed knot DGA has a similar
natural linearization, with respect to which we can deduce the
Alexander polynomial of the knot.

Let $K$ be a knot, and let $(\A,\d)$ denote the DGA over
$\Z[\mu^{\pm 1}]$ obtained by setting $\lambda=1$ in the framed knot
DGA of $K$. To define a linearization of $(\A,\d)$, we need an
augmentation of $(\A,\d)$ over $\Z[\mu^{\pm 1}]$, that is, an
algebra map $\epsilon:\thinspace \A\rightarrow \Z[\mu^{\pm 1}]$
which is $0$ in nonzero degrees and satisfies $\epsilon\circ\d = 0$.
In this case, an augmentation is just a ring homomorphism from the
cord algebra of $K$ with $\lambda=1$ to $\Z[\mu^{\pm 1}]$. There is
a natural choice for such a homomorphism, namely the one which sends
every cord to $1+\mu$; this clearly satisfies the relations for the
cord algebra.

If $\eps$ denotes the corresponding augmentation, which we call the
\textit{distinguished augmentation} of the framed knot DGA, then we
can construct a linearized version of $(\A,\d)$ in the standard way.
Let $\M$ be the subalgebra of $\A$ generated by all words (in
$a_{ij}$, $b_{ij}$, $c_{ij}$, $d_{ij}$, $e_i$) of length at least
$1$, and let $\varphi_\eps:\thinspace \A\rightarrow\A$ be the
algebra isomorphism sending each generator $g$ of $\A$ to
$g+\eps(g)$. The differential $\varphi_\eps \circ \d \circ
\varphi_\eps^{-1}$ maps $\M$ to itself and hence descends to a
differential $\d^{\lin}_{\eps}$ on the free $\Z[\mu^{\pm 1}]$-module
$\M/\M^2$. We define the \textit{linearized framed contact homology}
of $K$, written $HC_*^{\lin}(K)$, to be the graded homology of
$\d^{\lin}_{\eps}$.

To relate linearized homology to the Alexander polynomial, we recall
some elementary knot theory adapted from \cite{Rol}. Fix a knot $K$
in $S^3$, and let $\tilde{X}_K$ denote the infinite cyclic cover of
the knot complement. Then $H_1(\tilde{X}_K)$ has a natural
$\Z[t^{\pm 1}]$ module structure, with respect to which it is called
the \textit{Alexander invariant}. The Alexander polynomial
$\Delta_K(t)$ can be deduced from the Alexander invariant in the
usual way; in particular, if the Alexander invariant is cyclic of
the form $\Z[t^{\pm 1}]/(p(t))$, then $p(t)$ is the Alexander
polynomial up to units.

\begin{proposition} \label{prop:linAlex}
If $K$ is a knot, then we have an isomorphism of $\Q[\mu^{\pm
1}]$-modules
\[
HC_1^{\lin}(K)\otimes\Q \cong (H_1(\tilde{X}_K) \oplus \Q[\mu^{\pm
1}]) \otimes (H_1(\tilde{X}_K) \oplus \Q[\mu^{\pm 1}]) \oplus
(\Q[\mu^{\pm 1}])^m
\]
for some $m\geq 0$. Here $H_1(\tilde{X}_K)$ is understood to be a
$\Q[\mu^{\pm 1}]$-module by setting $t=-\mu$.
\end{proposition}

Before proving for Proposition~\ref{prop:linAlex}, we note that the
result, along with some simple linear algebra, implies that we can
calculate the Alexander polynomial from $HC_1^{\lin}$.

\begin{corollary} \label{cor:Alex}
The framed knot DGA with $\lambda=1$, along with the distinguished
augmentation of the DGA, determines the Alexander polynomial.
\end{corollary}

\begin{proof}
Suppose that the invariant factors for the $\Q[\mu^{\pm 1}]$-module
$H_1(\tilde{X}_K)$ are $p_1,\dots,p_k\in \Q[\mu^{\pm 1}]$, with
$p_1|\cdots|p_k$, so that $p_1\dots p_k$ is some rational multiple
of the Alexander polynomial. Then the invariant factors of
$(H_1(\tilde{X}_K) \oplus \Q[\mu^{\pm 1}]) \otimes (H_1(\tilde{X}_K)
\oplus \Q[\mu^{\pm 1}]) \oplus (\Q[\mu^{\pm 1}])^m$ are
$\gcd(p_i,p_j)$ for $1\leq i,j\leq k$, along with $2k$ additional
invariant factors which consist of $p_1,\dots,p_k$ with multiplicity
$2$ each. In other words, the invariant factors of $HC_1^{\lin}(K)$
are $p_1$ with multiplicity $2k+1$, $p_2$ with multiplicity $2k-1$,
\ldots, $p_k$ with multiplicity $3$. From $HC_1^{\lin}(K)$ we can
hence deduce $p_1,\dots,p_k$. This gives us the Alexander polynomial
up to a constant factor, which can be eliminated via the
normalization $\Delta_K(1)=1$.
\end{proof}

The proof of Proposition~\ref{prop:linAlex} will occupy the
remainder of this section. It is nearly identical to the proof of
the analogous Proposition~7.8 from \cite{Ng1} and is rather
unilluminating, and the reader may wish to skip it.

\begin{figure}
\centerline{
\includegraphics[height=0.7in]{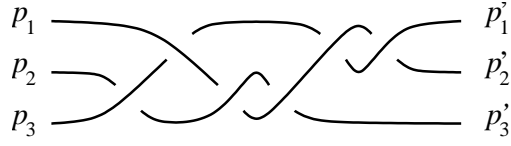}
}
\caption{
The braid $B=\sigma_1^2\sigma_2^2\sigma_1^{-1}\sigma_2$, with $v_B =
(1,\mu^2,1)$ and $u_B = (\mu^2,\mu^2,1)$.
}
\label{fig:braidex}
\end{figure}

We first give an explicit description of the distinguished
augmentation for a framed knot DGA, as derived by inspecting the
translation between cords and framed knot contact homology from the
proof of Theorem~\ref{thm:DGAcord}. Let $B\in B_n$ be a braid which
closes to a knot. Draw $B$ \textit{from right to left} with
$\sigma_k$ giving a crossing where strand $k$ crosses over strand
$k+1$, viewed from the right; see Figure~\ref{fig:braidex}. Write
$p_1,\ldots,p_n$ for the leftmost endpoints of the braid, and
$p_1',\ldots,p_n'$ for the rightmost endpoints. For $1\leq i\leq n$,
traverse the braid from left to right, beginning at $p_i$ and
looping around (that is, identify $p_j'$ with $p_j$) until $p_1$ is
reached; this traces out a path on the braid, which we define to be
trivial when $i=1$. Now let $r_i$ be the signed number of
undercrossings traversed by this path, where an undercrossing
contributes $1$ or $-1$ depending on whether it corresponds to a
braid generator $\sigma_k$ or $\sigma_k^{-1}$. Finally, define $v_B$
to be the vector $(\mu^{r_1},\ldots,\mu^{r_n})$.

Recall from Section~\ref{ssec:pfDGAinv} that $\Delta(v)$ is a
diagonal matrix whose diagonal entries are the entries of $v$, and
let $\Theta$ denote the $n\times n$ matrix whose entries are all
$-\mu-1$.

\begin{definition} \label{def:distaug}
The \textit{distinguished augmentation} of the framed knot DGA of a
braid $B\in B_n$ is the map which sends
\[
A \mapsto \Delta(v_B)\cdot\Theta\cdot\Delta(v_B)^{-1};
\]
in other words, it sends $a_{ij}$ to $\mu^{r_i-r_j}(-\mu-1)$ if
$i>j$ and $\mu^{r_i-r_j-1}(-\mu-1)$ if $i<j$.
\end{definition}

We will see shortly that the distinguished augmentation is, in fact,
an augmentation.

A bit more notation: let $\Lambda_B$ denote the matrix $\Lambda$
with $\lambda$ set to $\mu^{-w(B)}$ where $w(B)$ is the writhe of
$B$; that is, $\Lambda_B = \Delta((\mu^{-w(B)},1,\ldots,1))$. Define
the vector $u_B$ to be $(\mu^{s_1},\ldots,\mu^{s_n})$, where $s_i$
is the signed number of undercrossings traversed if we begin at
$p_i$, travel from left to right, and end when we reach the right
hand end of the braid; see Figure~\ref{fig:braidex}.

Also, let $\hat{B}$ be the braid obtained by reversing the word
which gives $B$. Finally, let $\Bur{}$ denote the Burau
representation on $B_n$, defined over $\Z[t^{\pm 1}]$, so that
$\Bur{\sigma_k}$ is the $n\times n$ matrix which is the identity
except from the $2\times 2$ submatrix formed by the $k,k+1$ rows and
columns, which is
$\left( \begin{smallmatrix} 1-t & t \\ 1 & 0 \end{smallmatrix}
\right)$.
We now give a series of lemmas, culminating in the proof
of Proposition~\ref{prop:linAlex}.

\begin{lemma}
Let $v_k$ denote the vector $(1,\ldots,1,\mu,1,\ldots,1)$, with
$\mu$ in the $k$-th position. Then
\begin{align*}
\hom{\sigma_k}(\Theta) &= \Delta(v_k) \cdot \Theta_k \cdot
\Delta(v_k)^{-1},\\
\Phil{\sigma_k}(\Theta) &= \Delta(v_k) \cdot \left(
\Bur{\sigma_k}|_{t=-1/\mu} \right),\\
\Phir{\sigma_k}(\Theta) &= \left( \Bur{\sigma_k}|_{t=-\mu} \right)^T
\cdot \Delta(v_k)^{-1},
\end{align*}
where $^T$ denotes transpose.
\label{lem:aug1}
\end{lemma}

\begin{proof}
Easy calculation using the definitions, along with the expressions
for $\Phil{\sigma_k}$ and $\Phir{\sigma_k}$ from \cite[Lemma
4.6]{Ng1}.
\end{proof}

\begin{lemma}
We have
\[
\Phil{B}(\Theta) = \Delta(u_B) \cdot \left( \Bur{\hat{B}}
|_{t=-1/\mu} \right) = \Lambda_B^{-1} \cdot \Delta(s(B) v_B)^{-1}
\cdot \Delta(v_B) \cdot \left( \Bur{\hat{B}} |_{t=-1/\mu} \right)
\]
and
\[
\Phir{B}(\Theta) = \left( \Bur{\hat{B}}|_{t=-\mu} \right) \cdot
\Delta(u_B)^{-1} = \left( \Bur{\hat{B}}|_{t=-\mu} \right) \cdot
\Delta(v_B)^{-1} \cdot \Delta(s(B)v_B) \cdot \Lambda_B,
\]
with $s(B)$ defined as in Section~\ref{ssec:pfDGAinv}.
\label{lem:aug2}
\end{lemma}

\begin{proof}
We will establish the formula for $\Phil{B}(\Theta)$, with
$\Phir{B}(\Theta)$ established similarly. The second equality
follows from the fact that $s_i = r_i - r_{s(B)(i)}$ for $i\geq 2$
and $s_1 = r_1 - r_{s(B)(1)} + w(B)$. To prove the first equality,
we use induction, noting that $u_B$ can be defined even when $B$
closes to a multi-component link rather than a knot.

By \cite[Prop. 4.4]{Ng1}, Lemma~\ref{lem:lambda}, and
Lemma~\ref{lem:aug1}, we have
\begin{align*}
\Phil{\sigma_k B}(\Theta) &= \Phil{B}(\Delta(v_k) \cdot \Theta \cdot
\Delta(v_k)^{-1}) \cdot \Phil{\sigma_k}(\Theta) \\
&= \Delta(s(B)v_k) \cdot \Phil{B}(\Theta) \cdot \Delta(v_k)^{-1} \\
&= \Delta(u_{\sigma_k B}) \cdot \Delta(u_B)^{-1} \cdot
\Phil{B}(\Theta) \cdot \left( \Bur{\sigma_k}|_{t=-1/\mu} \right).
\end{align*}
Hence the formula holds for $B$ if and only if it holds for
$\sigma_k B$. Since it certainly holds when $B$ is the trivial
braid, the formula holds for general $B$ by induction.
\end{proof}

\begin{lemma}
The distinguished augmentation is in fact an augmentation for the
framed knot DGA of $B$, once we set $\lambda=\mu^{-w(B)}$.
\label{lem:aug3}
\end{lemma}

\begin{proof}
Define $\Theta_B = \Delta(v_B) \cdot \Theta \cdot \Delta(v_B)^{-1}$.
We want to show that the map $A \mapsto \Theta_B$ sends the
differentials of the degree $1$ generators of the framed knot DGA to
$0$. By Lemmas~\ref{lem:lambda} and \ref{lem:aug2}, we have
\[
\Phil{B}(\Theta_B) = \Delta(s(B)v_B) \cdot \Phil{B}(\Theta) \cdot
\Delta(v_B)^{-1} = \Lambda_B^{-1} \cdot \Delta(v_B) \cdot \left(
\Bur{\hat{B}}|_{t=-1/\mu} \right) \cdot \Delta(v_B)^{-1}
\]
and thus
\[
(1-\Lambda_B \cdot \Phil{B}(\Theta_B))\cdot \Theta_B = \Delta(v_B)
\cdot \left(1- \Bur{\hat{B}}|_{t=-1/\mu} \right) \cdot \Theta \cdot
\Delta(v_B)^{-1} = 0;
\]
similarly $\Theta_B \cdot (1-\Phir{B}(\Theta_B)\cdot \Lambda_B^{-1})
= 0$.
\end{proof}

\begin{proof}[Proof of Proposition~\ref{prop:linAlex}]
Let $B$ be a braid whose closure is $K$, and consider the
linearization of the modified framed knot DGA of $B$, as defined in
Section~\ref{ssec:pfDGAinv}. By the calculation from the proof of
Lemma~\ref{lem:aug3}, we have
\[
\d^{\lin}_{\eps} F = B - C - \Lambda_B \cdot \Phil{B}(\Theta_B)\cdot
D = B - C - \Delta(v_B) \cdot \left(\Bur{\hat{B}}|_{t=-1/\mu}\right)
\cdot \Delta(v_B)^{-1} \cdot D
\]
and similarly $\d^{\lin}_{\eps} E = B - D - C\cdot \Delta(v_B)\cdot
\left(\Bur{\hat{B}}|_{t=-\mu}\right) \cdot \Delta(v_B)^{-1}$. If we
change variables by replacing $B$ by $\Delta(v_B)\cdot
B\cdot\Delta(v_B)^{-1}$ and similarly for $C,D,E,F$, then the
differential becomes $\d^{\lin}_{\eps} E = B - D - C\cdot
\left(\Bur{\hat{B}}|_{t=-1/\mu}\right)$, $\d^{\lin}_{\eps} F = B - C
-\left(\Bur{\hat{B}}|_{t=-\mu}\right)\cdot D$.

The matrices $\Bur{\hat{B}}|_{t=-1/\mu}$ and
$\Bur{\hat{B}}|_{t=-\mu}$ both give presentations over $\Q[t^{\pm
1}]$ for the Alexander invariant of $K$, or, more precisely, for the
direct sum of the Alexander polynomial and a free summand $\Q[t^{\pm
1}]$, once we identify $t$ with $-\mu$. (Note that the Alexander
invariant is unchanged if we invert $t$.) It is now straightforward
to deduce the proposition; just use the method of the proof of
\cite[Prop.~7.8]{Ng1}, working over the ring $\Q[\mu^{\pm 1}]$
rather than $\Z$.
\end{proof}

\secdivide
%*********************************************************************
%*********************************************************************
\section{Invariants from Augmentation}
\label{sec:aug}

In this section, we focus on applications of framed knot contact
homology. The essential difficulty is that it is a nontrivial task
to determine when two DGAs over $\ring$ are equivalent. To obtain
some computable invariants, we consider augmentations of the DGA.
These give rise to a large family of readily calculable numerical
invariants of knots (Section~\ref{ssec:augno}), as well as a
$2$-variable polynomial with very strong ties to the $A$-polynomial
(Section~\ref{ssec:augpoly}).

\subsecdivide
%*********************************************************************
\subsection{Augmentation numbers}
\label{ssec:augno}

Let $\k$ be a field, and $\k^* = \k\setminus\{0\}$. Any pair
$(\lambda_0,\mu_0)\in (\k^*)^2$ gives rise to a map $\ring \otimes
\k \rightarrow \k$ by sending $\lambda$ to $\lambda_0$ and $\mu$ to
$\mu_0$. If $(\A,\d)$ is a DGA over $\ring$, then
$(\lambda_0,\mu_0)$ allows us to project $\A\otimes\k$, which is an
algebra over $\k[\lambda^{\pm 1},\mu^{\pm 1}]$, to an algebra over
$\k$, which we write as $\A|_{(\lambda_0,\mu_0)}$.

\begin{definition} \label{def:aug}
An \textit{augmentation} of a DGA $(\A,\d)$ over a field $\k$ is an
algebra map $\vareps:\thinspace \A\rightarrow \k$ such that
$\vareps(1)=1$, $\vareps(a)=0$ if $a\in\A$ has pure nonzero degree,
and $\vareps\circ\d = 0$.
\end{definition}

If $\k$ is finite and $\A$ is finitely generated as an algebra, then
there are only finitely many possible augmentations.

\begin{definition} \label{def:augno}
Let $K$ be a knot with framed DGA $(\A,\d)$, let $\k$ be a finite
field, and let $(\lambda_0,\mu_0)\in (\k^*)^2$. The
\textit{augmentation number} $\Aug(K,\k,\lambda_0,\mu_0)$ is the
number of augmentations of $\A|_{(\lambda_0,\mu_0)}$ over $\k$.
\end{definition}

\noindent Since a knot has many different framed DGAs depending on
its braid representative, Definition~\ref{def:augno} technically
requires proof of invariance; however, since framed knot DGAs have
generators only in nonnegative degree, it is easy to see that an
augmentation of $\A|_{(\lambda_0,\mu_0)}$ is the same as an algebra
map $(HC_0(K)\otimes\k)|_{\lambda=\lambda_0,\mu=\mu_0} \rightarrow
\k$, and so augmentation numbers are well-defined.

From Section~\ref{ssec:cordDGA}, if the bridge number of a knot $K$
is $m$, then there is a way to write $HC_0(K)$ as an algebra with
$m(m-1)$ generators and $2m^2$ relations. Calculating an
augmentation number $\Aug(K,\k,\lambda_0,\mu_0)$ then involves
counting how many of the possible $|\k|^{m(m-1)}$ maps of the
generators to $\k$ satisfy the relations; this can be readily done
by computer when $|\k|$ is small enough.

In particular, if $p$ is a prime, the finite field $\Z_p$ gives rise
to $(p-1)^2$ augmentation numbers of a knot, corresponding to the
different choices of $(\lambda_0,\mu_0)$, each of which is a knot
invariant. It can be calculated without too much effort that nearly
any pair of nonisotopic knots with $8$ or fewer crossings can be
distinguished by one of the $69$ total augmentation numbers for
$\Z_p$ with $p=2,3,5,7$. The only pairs not so distinguished are the
knot $8_{15}$ and its mirror, and the connect sum $3_1\# 5_1$ and
its mirror; each of these pairs share the same augmentation numbers
over fields up to $\Z_7$ but can be distinguished using augmentation
numbers over $\Z_{11}$.

\begin{figure}
\centerline{ \includegraphics[width=5in]{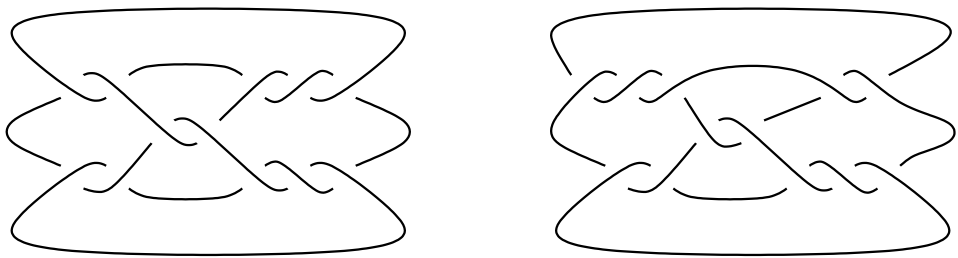} }
\caption{The Kinoshita--Terasaka knot (left) and Conway mutant
(right).} \label{fig:KTConway}
\end{figure}

As an application of augmentation numbers, we can calculate that the
Kinoshita--Terasaka knot and its Conway mutant, as shown in
Figure~\ref{fig:KTConway}, do not have all the same augmentation
numbers: $\Aug(\textrm{KT},\Z_7,3,5) = 2$ while
$\Aug(\textrm{Conway},\Z_7,3,5) = 1$. The author's web page includes
a \textit{Mathematica} program to calculate augmentation numbers.

\begin{proposition} \label{prop:KTconway}
The Kinoshita--Terasaka knot and its Conway mutant have
nonisomorphic cord algebras (and thus inequivalent framed knot
DGAs).
\end{proposition}

We conclude that augmentation numbers are a rather effective tool to
distinguish between knots. It is not known to the author whether
there are two nonisotopic knots which share the same augmentation
numbers for all fields $\Z_p$.

\subsecdivide
%*********************************************************************
\subsection{The augmentation polynomial}
\label{ssec:augpoly}

Here we consider augmentations, as in the previous section, but now
over the field $\C$. If we fix $(\lambda_0,\mu_0)\in(\C^*)^2$ and a
knot $K$, then we obtain a DGA $(\A|_{(\lambda_0,\mu_0)},\d)$ over
$\C$, and we can ask whether this DGA has an augmentation, i.e.,
whether there is an algebra map $(HC_0(K) \otimes
\C)|_{\lambda=\lambda_0,\mu=\mu_0} \rightarrow \C$. Since $HC_0(K)$
is finitely generated and finitely presented, this is just a
question in elimination theory, and the (closure of the) locus of
$(\lambda_0,\mu_0)$ for which the DGA has an augmentation is an
algebraic set.

\begin{definition}
Let $K$ be a knot, with framed knot DGA $(\A,\d)$. The set
$\{(\lambda_0,\mu_0)\in(\C^*)^2:\thinspace
(\A|_{(\lambda_0,\mu_0)},\d) \textrm{ has an augmentation}\}$ is
called the \textit{augmentation variety} of $K$. If the augmentation
variety is not $2$-dimensional (and thus a Zariski open subset of
$(\C^*)^2$), then the union of its $1$-dimensional components is the
zero set of the \textit{augmentation polynomial}
$\tilde{A}_K(\lambda,\mu)\in\C[\lambda,\mu]$, which is unique up to
constant multiplication once we specify that it has no repeated
factors and is not divisible by $\lambda$ or $\mu$.
\end{definition}

\noindent The author does not know of any examples for which the
augmentation variety has any components which are not
$1$-dimensional.

\begin{conjecture} \label{conj:augpoly}
The augmentation polynomial is defined for all knots in $S^3$.
\end{conjecture}

\begin{proposition} \label{prop:augpoly1}
The augmentation polynomial $\tilde{A}_K(\lambda,\mu)$ is always
divisible by $(\lambda-1)(\mu+1)$.
\end{proposition}

\begin{proof}
If $\lambda_0=1$, then there is a map from $HC_0(K)|_{(1,\mu_0)}$ to
$\C$ sending all cords to $1+\mu$; if $\mu_0=-1$, then there is a
trivial map from $HC_0(K)|_{(\lambda_0,-1)}$ to $\C$ sending all
cords to $0$.
\end{proof}

We computed in Section~\ref{ssec:htpydef} that the cord algebra for
the unknot is $\ring/((\lambda-1)(\mu+1))$; hence the unknot has
augmentation polynomial
\[
\tilde{A}_0(\lambda,\mu) = (\lambda-1)(\mu+1).
\]
A slightly more interesting case is the right hand trefoil, for
which we computed in Section~\ref{ssec:cordDGA} that the cord
algebra is $(\ring)[x]/(x^2-x-\lambda\mu^2-\lambda\mu,x^2-\lambda\mu
x-\lambda-\lambda\mu)$. For this to have a map to $\C$, we need the
two polynomials generating the ideal to have a common root. Their
resultant is $\lambda(\lambda-1)(\mu+1)(1-\lambda\mu^3)$, and so the
right hand trefoil has augmentation polynomial
\[
\tilde{A}_{3_1}(\lambda,\mu) = (\lambda-1)(\mu+1)(1-\lambda\mu^3).
\]

We will give some more computations of augmentation polynomials
shortly, but first we establish a couple more properties of
augmentation polynomials. For Laurent polynomials in $\ring$, let
$\doteq$ denote equality up to multiplication by a unit in $\ring$.

\begin{proposition} \label{prop:augpoly2}
The augmentation polynomial is independent of the orientation of the
knot. If $K$ is a knot with mirror $\overline{K}$, then
\[
\tilde{A}_K(\lambda,\mu) \doteq \tilde{A}_K(\lambda^{-1},\mu^{-1})
\hspace{4ex} \text{and} \hspace{4ex}
\tilde{A}_{\overline{K}}(\lambda,\mu) \doteq
\tilde{A}_K(\lambda,\mu^{-1}).
\]
\end{proposition}

\begin{proof}
Follows immediately from Proposition~\ref{prop:cordalgproperties};
note that the augmentation polynomial depends only on the
commutative framed cord algebra, rather than the noncommutative
version.
\end{proof}

\noindent Note also that one can define an augmentation polynomial
for knots with arbitrary framing $f$ by using the $f$-framed knot
DGA instead of the $0$-framed knot DGA; then the augmentation
polynomial of a knot $K$ with framing $f$ is
$\tilde{A}_K(\lambda\mu^{-f},\mu)$.

\begin{proposition} \label{prop:augpolyconnectsum}
The augmentation variety $V_{K_1\# K_2}$ for a connect sum $K_1\#
K_2$ is related to the augmentation varieties $V_{K_1},V_{K_2}$ for
$K_1,K_2$ as follows:
\[
V_{K_1\# K_2} = \{(\lambda_1\lambda_2,\mu)\in (\C^*)^2\,|\,
(\lambda_1,\mu)\in V_{K_1} \text{ and } (\lambda_2,\mu)\in
V_{K_2}\}.
\]
This formula defines $\tilde{A}_{K_1\# K_2}$ in terms of
$\tilde{A}_{K_1}$ and $\tilde{A}_{K_2}$. In particular, both
$\tilde{A}_{K_1}$ and $\tilde{A}_{K_2}$ divide $\tilde{A}_{K_1\#
K_2}$.
\end{proposition}

\begin{proof}
Since $\{\mu=-1\}$ is a component of all augmentation varieties, it
suffices to prove the equality in the statement in the complement of
this component.

By Van Kampen's Theorem, $\pi_1(S^3\setminus (K_1\# K_2))$ is the
free product of $\pi_1(S^3\setminus K_1)$ and $\pi_1(S^3\setminus
K_2)$, modulo the identification of the two meridians. Thus there
are elements $l_1,l_2,m \in \pi_1(S^3\setminus (K_1\# K_2))$ such
that $l_1l_2,m$ are the longitude and meridian classes of $K_1\#
K_2$, and $l_i,m$ project to the longitude and meridian classes in
$\pi_1(S^3\setminus K_i)$ for $i=1,2$.

A map $p:\thinspace HC_0(K_1\#
K_2)|_{(\lambda_0,\mu_0)}\rightarrow\C$ with $\mu_0\neq -1$ sends
$[l_1],[l_2]$ to complex numbers which we write as
$\lambda_{0,1}(1+\mu_0),\lambda_{0,2}(1+\mu_0)$. We claim that $p$
induces maps from $HC_0(K_i)|_{(\lambda_{0,i},\mu_0)}$ to $\C$. The
only relation for $HC_0(K_i)|_{(\lambda_{0,i},\mu_0)}$ from
Definition~\ref{def:knotcordalg} that is not trivially preserved by
$p$ is $[\gamma l_i] = [l_i \gamma] = \lambda_{0,i} [\gamma]$. But
in $HC_0(K_1\# K_2)|_{(\lambda_0,\mu_0)}$, we have
\[
[\gamma] [l_i] = [\gamma l_i]+[\gamma m l_i] = [\gamma l_i] +
[\gamma l_i m] = (1+\mu_0) [\gamma l_i],
\]
and so $p([\gamma l_i]) = p([\gamma] [l_i])/(1+\mu_0) =
\lambda_{0,i} p([\gamma])$. Similarly $p([l_i \gamma]) =
\lambda_{0,i} p([\gamma])$, and so $p$ does indeed restrict to maps
from $HC_0(K_i)|_{(\lambda_{0,i},\mu_0)}$ to $\C$. In addition, we
have
\[
\lambda_0(1+\mu_0) = p([l]) = p([l_1l_2]) = \lambda_{0,1} p([l_2]) =
\lambda_{0,1}\lambda_{0,2}(1+\mu_0).
\]
We conclude that if $(\lambda_0,\mu_0) \in V_{K_1\# K_2}$, then
there exist $\lambda_{0,1},\lambda_{0,2}$ with
$\lambda_{0,1}\lambda_{0,2}=\lambda_0$ such that
$(\lambda_{0,i},\mu_0)\in V_{K_i}$ for $i=1,2$.

Conversely, suppose that we are given maps $p_i:\thinspace
HC_0(K_i)|_{(\lambda_{0,i},\mu_0)} \rightarrow \C$. We construct a
map $p:\thinspace HC_0(K_1\#
K_2)|_{(\lambda_{0,1}\lambda_{0,2},\mu_0)} \rightarrow \C$ as
follows: if $\gamma_{j,i}\in\pi_1(S^3\setminus K_i)$ for $i=1,2$ and
$1\leq j\leq k$, then
\[
p([\gamma_{1,1}\gamma_{1,2}\cdots\gamma_{k,1}\gamma_{k,2}]) =
\frac{p_1([\gamma_{1,1}])p_2([\gamma_{1,2}])\cdots
p_1([\gamma_{k,1}])p_2([\gamma_{k,2}])}{(1+\mu_0)^{k-1}}.
\]
Since any element of $\pi_1(S^3\setminus (K_1\# K_2))$ can be
written in the form
$\gamma_{1,1}\gamma_{1,2}\cdots\gamma_{k,1}\gamma_{k,2}$, this
defines $p$ on $HC_0(K_1\#
K_2)|_{(\lambda_{0,1}\lambda_{0,2},\mu_0)}$. It is straightforward
to check that $p$ preserves the relations defining $HC_0(K_1\#
K_2)|_{(\lambda_{0,1}\lambda_{0,2},\mu_0)}$ and is hence
well-defined.
\end{proof}

To illustrate Proposition~\ref{prop:augpolyconnectsum}, we can
calculate that the augmentation polynomial for the connect sum of
two right hand trefoils is
\[
\tilde{A}_{3_1\# 3_1}(\lambda,\mu) =
(\lambda-1)(\mu+1)(1-\lambda\mu^3)(1-\lambda\mu^6),
\]
while the augmentation polynomial for the connect sum of a left hand
trefoil and a right hand trefoil is
\[
\tilde{A}_{3_1\#\overline{3_1}}(\lambda,\mu) =
(\lambda-1)(\mu+1)(1-\lambda\mu^3)(\lambda-\mu^3).
\]

\subsecdivide
%*********************************************************************
\subsection{The augmentation polynomial and the $A$-polynomial}
\label{ssec:augpolyApoly}

The augmentation polynomial is very closely related to the
$A$-polynomial introduced by Cooper et al.\ in \cite{CCGLS}, which
we briefly review here. Let $K\in S^3$ be a knot, and let $l,m$
denote the longitude and meridian classes in $\pi_1(S^3\setminus
K)$, as usual. Consider a representation $\rho:\thinspace
\pi_1(S^3\setminus K) \rightarrow SL_2\C$. Since $l,m$ commute,
$\rho$ is conjugate to a representation in which $\rho(l)$ and
$\rho(m)$ are upper triangular, and so we restrict our attention to
representations $\rho$ with this property. There is a map from
$SL_2\C$-representations of $\pi_1(S^3\setminus K)$ to $(\C^*)^2$
sending $\rho$ to $((\rho(l))_{11},(\rho(m))_{11})$, where $M_{11}$
denotes the top left entry of $M$, and the $1$-dimensional part of
the image of this map forms a variety, which is the zero set of the
$A$-polynomial $A_K\in\C[\lambda,\mu]$. (Note that we are thus
assuming the convention that $\lambda-1$ always divides $A_K$.)

The following result shows that the $A$-polynomial is subsumed in
the $\tilde{A}$-polynomial. Note that the $A$-polynomial involves
only even powers of $\mu$ \cite[Prop.\ 2.9]{CCGLS}.

\begin{proposition} \label{prop:Apoly}
Let $A_K(\lambda,\mu)$, $\tilde{A}_K(\lambda,\mu)$ denote the
$A$-polynomial and augmentation polynomial of a knot $K$. Then
\[
(1-\mu^2) A_K(\lambda,\mu)\,\left|\, \tilde{A}_K(\lambda,-\mu^2).
\right.
\]
\end{proposition}

\begin{proof}
From \cite[\S 2.8]{CCGLS}, $\mu\pm 1$ does not divide the
$A$-polynomial. On the other hand, by
Proposition~\ref{prop:augpoly1}, $\tilde{A}_K(\lambda,-\mu^2)$ is
divisible by $1-\mu^2$. It thus suffices to prove that
$A_K(\lambda,\mu)$ divides $\tilde{A}_K(\lambda,-\mu^2)$.

Suppose that we have an $SL_2\C$-representation $\rho$ of
$\pi_1(S^3\setminus K)$ which maps to $(\lambda_0,\mu_0)$ in
$(\C^*)^2$. We wish to show that there is a map from the cord
algebra of $K$ with $(\lambda,\mu)=(\lambda_0,-\mu_0^2)$ to $\C$.
Since $\mu+1$ is a factor of $\tilde{A}_K(\lambda,\mu)$, we may
assume that $\mu_0 \neq \pm 1$. Then we can conjugate $\rho$ so that
it becomes diagonal on the peripheral subgroup:
\[
\rho(l) = \left( \begin{matrix} \lambda_0 & 0 \\ 0 & \lambda_0^{-1}
\end{matrix} \right)
\hspace{0.5in} \text{and} \hspace{0.5in} \rho(m) = \left(
\begin{matrix} \mu_0 & 0 \\ 0 & \mu_0^{-1}
\end{matrix} \right).
\]

We use the formulation of the cord algebra from
Definition~\ref{def:knotcordalg}, with $\lambda$ replaced by
$\lambda_0$ and $\mu$ replaced by $-\mu_0^2$; we can change
variables by replacing each generator $[\gamma]$ by
$(1-\mu_0^2)(-\mu_0)^{\lk(\gamma,K)}[\gamma]$ to obtain an algebra
with defining relations
\begin{enumerate}
\item $[e]=1$;
\item $[\gamma m]=[m\gamma] = \mu_0[\gamma]$ and
$[\gamma l]=[l \gamma] = \lambda_0[\gamma]$;
\item $[\gamma_1\gamma_2] - \mu_0 [\gamma_1 m \gamma_2] = (1-\mu_0^2)
[\gamma_1] [\gamma_2]$.
\end{enumerate}
Now if we send each generator $[\gamma]$ to the top left entry of
$\rho(\gamma)$, it is easy to check that the above relations are
preserved, and so we obtain a map from the cord algebra to $\C$ with
$(\lambda,\mu)=(\lambda_0,-\mu_0^2)$, as desired.
\end{proof}

Proposition~\ref{prop:Apoly} has some immediate consequences. The
first uses a result of Dunfield and Garoufalidis \cite{DG} which is
based on work of Kronheimer and Mrowka \cite{KM} on knot surgery and
$SU_2$ representations.

\begin{proposition}
\label{prop:nontrivial} The cord algebra distinguishes
the unknot from any other knot in $S^3$.
\end{proposition}

\begin{proof}
We have seen that the cord algebra for the unknot is
$\ring/((\lambda-1)(\mu+1))$, and so the augmentation variety for
the unknot is the union of the lines $\lambda=1$ and $\mu=-1$. On
the other hand, the main result of \cite{DG} states that any
nontrivial knot in $S^3$ has nontrivial $A$-polynomial (i.e., not
equal to $\lambda-1$). Hence the augmentation variety for a
nontrivial knot either is $2$-dimensional or contains a nontrivial
$1$-dimensional component; in either case, it strictly contains the
augmentation variety for the unknot.
\end{proof}

The next consequence of Proposition~\ref{prop:Apoly} relates the
augmentation variety to the Alexander polynomial. By analogy with
\cite{CCGLS}, say that a point $(\lambda_0,\mu_0)$ in the
augmentation variety of $K$ is \textit{reducible} if $\lambda_0=1$
and the corresponding map from the cord algebra (from
Definition~\ref{def:cordcordalg}) to $\C$ sends every cord to
$1+\mu_0$. (This is slightly bad notation: we mean a point which is
on a component of the augmentation variety different from
$\lambda-1$, with the corresponding map from the cord algebra to
$\C$.)

\begin{proposition}
If $-\mu_0$ is a root of the Alexander polynomial $\Delta_K(t)$,
then $(1,\mu_0)$ is a reducible point on a nontrivial component of
the augmentation variety of $K$.
\end{proposition}

\begin{proof}
This follows directly from \cite[Prop.\ 6.2]{CCGLS}.
\end{proof}

We now examine the difference between the augmentation polynomial
and the $A$-polynomial for some specific knots. It seems that the
two polynomials essentially coincide for $2$-bridge knots, but as of
this writing, a proof has not been completed, and so the following
``result'' is presented as a conjecture.

\begin{conjecture}
\label{conj:Apoly2bridge} If $K$ is a $2$-bridge knot, then
\[
(1-\mu^2)A_K(\lambda,\mu) \doteq \tilde{A}_K(\lambda,-\mu^2).
\]
\end{conjecture}

%\begin{proof}
%\end{proof}

On the other hand, there are knots for which
$\tilde{A}_K(\lambda,-\mu^2)$ contains nontrivial factors not
present in $A_K(\lambda,\mu)$. Define the polynomial
$B_K(\lambda,\mu)$, defined up to multiplication by constants, to be
the product of the factors of $\tilde{A}_K(\lambda,\mu)$, besides
$\mu+1$, which do not have a corresponding factor in
$A_K(\lambda,\mu)$; that is,
\[
B_K(\lambda,-\mu^2) =
\frac{\tilde{A}_K(\lambda,-\mu^2)}{(1-\mu^2)A_K(\lambda,\mu)}.
\]
Conjecture~\ref{conj:Apoly2bridge} states that $B_K = 1$ when $K$ is
a $2$-bridge knot.

Suppose that $K$ is the torus knot $T(3,4)$. The group
$\pi_1(S^3\setminus T(3,4))$ has presentation $\langle
x,y\,:\,x^3=y^4\rangle$, with peripheral classes $m=xy^{-1}$ and
$l=x^3m^{-12}$. A quick calculation then shows that
\[
A_{T(3,4)}(\lambda,\mu) =
(\lambda-1)(1+\lambda\mu^{12})(1-\lambda\mu^{12}).
\]
On the other hand, we can define the cord algebra for $T(3,4)$ from
Definition~\ref{def:knotDGAbraid}, using the braid
$(\sigma_1\sigma_2)^{-4}\in B_3$ which closes to $T(3,4)$. It is
straightforward to check that if we set $\lambda=\mu^{-8}$ and
$a_{12}=a_{21}=a_{13}=a_{31}=a_{23}=a_{32}=-1$, then the relations
defining the cord algebra for $T(3,4)$ vanish, and so
$1-\lambda\mu^8$ is a factor of $\tilde{A}_{T(3,4)}$. In fact, a
slightly more involved calculation using resultants demonstrates
that
\[
B_{T(3,4)}(\lambda,\mu) = 1-\lambda\mu^8 \not\doteq 1.
\]

Using \textit{Mathematica}, one can compute the augmentation
polynomials for many other knots. In particular, for the
non-two-bridge knots whose $A$-polynomials are computed in
\cite{CCGLS}, we have
\begin{align*}
B_{8_5}(\lambda,\mu) &= 1- \lambda\mu^4; \\
B_{8_{20}}(\lambda,\mu) &= 1- \lambda\mu^2; \\
B_{P(-2,3,7)}(\lambda,\mu) &= 1- \lambda\mu^{12}.
\end{align*}
The author does not know whether there is any geometric significance
to the $B$-polynomial, as there is for the $A$-polynomial.

\secdivide
%*********************************************************************
%*********************************************************************
\section{Extensions of the Invariant}
\label{sec:extensions}

In Section~\ref{ssec:htpydef}, the cord algebra was defined for any
connected cooriented codimension $2$ submanifold of any manifold. In
fact, the cord algebra, and sometimes the full DGA invariant, can be
defined for many geometric objects. Here we discuss extending the
invariant to four settings other than knots in $S^3$. What follows
is expository and somewhat sketchy at points.

\subsecdivide
%*********************************************************************
\subsection{Links in $S^3$}
\label{ssec:links}

The cord algebra can be defined for (oriented, framed) links in
$S^3$ as for knots. For an $n$-component link $L$, the base ring is
$\Z[\lambda_1^{\pm 1},\mu_1^{\pm 1},\ldots,\lambda_n^{\pm
1},\mu_n^{\pm 1}]$ since there are now $n$ longitudes and $n$
corresponding meridians in $\pi_1(S^3\setminus L)$, and the
relations in Definition~\ref{def:knotcordalg} change accordingly.
The cord algebra can be written in terms of cords and skein
relations as for knots:
\begin{gather*}
\raisebox{-0.17in}{\includegraphics[width=0.4in]{skein21.eps}}
= 1+\mu_i, \\
\raisebox{-0.17in}{\includegraphics[width=0.4in]{skein31.eps}} =
\lambda_i
\raisebox{-0.17in}{\includegraphics[width=0.4in]{skein32.eps}}
\hspace{0.4in} \textrm{and} \hspace{0.4in}
\raisebox{-0.17in}{\includegraphics[width=0.4in]{skein32-mod.eps}} =
\lambda_i
\raisebox{-0.17in}{\includegraphics[width=0.4in]{skein31-mod.eps}},
\\
\raisebox{-0.17in}{\includegraphics[width=0.4in]{skein12.eps}} +
\mu_i \raisebox{-0.17in}{\includegraphics[width=0.4in]{skein11.eps}}
= \raisebox{-0.17in}{\includegraphics[width=0.4in]{skein13.eps}}
\cdot
\raisebox{-0.17in}{\includegraphics[width=0.4in]{skein14.eps}},
\end{gather*}
where $i$ is the label of the link component depicted in the skein
relation.

It is possible to extend the framed knot DGA from
Section~\ref{ssec:cordDGA} to an invariant of links in such a way
that the degree $0$ homology of the DGA is the cord algebra of the
link; one just needs to supply $\lambda,\mu$ in the definition of
the DGA with indices corresponding to components of the link. There
is then a canonical augmentation of the DGA which sends each
$a_{ij}$ to $1+\mu_{n(i)}$, where $n(i)$ is the label of the link
component containing strand $i$. It should be the case that the
linearized homology of the DGA with respect to this augmentation
contains the multivariable Alexander polynomial of the link, just as
the linearized homology of the framed knot DGA contains the usual
Alexander polynomial of a knot; however, the author has not worked
through the detailed computation in general.

\subsecdivide
%*********************************************************************
\subsection{Spatial graphs}
\label{ssec:graphs}

As noted in \cite{Ng2}, the cord algebra over $\Z$ extends naturally
to an invariant of spatial graphs (embedded graphs in $S^3$) modulo
neighborhood equivalence; two spatial graphs are
\textit{neighborhood equivalent} (terminology by Suzuki) if they
have ambient isotopic tubular neighborhoods. If $\Gamma$ is a
spatial graph and $H_1(\Gamma)$ has rank $g$, then one can lift the
base ring of the cord algebra of $\Gamma$ from $\Z$ to
$\Z[\lambda_1^{\pm 1},\dots,\lambda_g^{\pm 1},\mu_1^{\pm
1},\dots,\mu_g^{\pm 1}]$, where the $\lambda_i$ represent generators
of $H_1(\Gamma)$ and the $\mu_i$ represent $g$ linearly independent
meridians. However, since the choices of longitudes and meridians
are not canonical in general, an isomorphism of cord algebras of
neighborhood equivalent graphs may not fix the base ring.

One can also generalize the framed knot DGA to spatial graphs by
using the formulation from Section~\ref{ssec:cordDGA}. Note that a
diagram of $\Gamma$ with, say, $n$ crossings will have $n-g+1$
connected components. The DGA for $\Gamma$ has $(n-g+1)(n-g)$
generators $a_{ij}$ in degree $0$, $2n(n-g+1)$ generators $b_{\alpha
i},c_{i\alpha}$ in degree $1$, and $n^2+n$ generators
$d_{\alpha\beta},e_\alpha$ in degree $2$; the differential is as in
Section~\ref{ssec:cordDGA}, where $\Psi^L,\Psi^L_1,\Psi^L_2$ are now
$n\times(n-g+1)$ matrices and $\Psi^R,\Psi^R_1,\Psi^R_2$ are
$(n-g+1)\times n$ matrices. The proof that the DGA is invariant
under neighborhood equivalence is identical to the corresponding
(omitted) proof of invariance for knots under isotopy.

By the definition of neighborhood equivalence, the cord algebra and
framed knot DGA also yield isotopy invariants for any closed
embedded surface in $S^3$ whose complement has a component which is
a handlebody; simply consider the spatial graph given by the core of
the handlebody. In general, any embedded surface in $\R^3$ (or
immersed surface without dangerous self-tangencies) should yield a
DGA and ``cord algebra'' from symplectic considerations, given by
the Legendrian contact homology of its conormal bundle (Legendrian
lift) in $ST^*\R^3$, but a combinatorial form for these invariants
is currently unknown. When the surface arises as the boundary of a
tubular neighborhood of a spatial graph, we conjecture that the DGA
for the spatial graph defined above calculates the Legendrian
contact homology of the lift of the associated surface.

\subsecdivide
%*********************************************************************
\subsection{Virtual and welded knots in $S^3$}
\label{ssec:virtual}

The framed knot DGA extends easily to virtual knots (for
definitions, see e.g. \cite{FKM}), via the formulation in
Section~\ref{ssec:cordDGA}. The definition of the invariant from
Section~\ref{ssec:cordDGA} carries over wholesale, with the
clarification that only actual crossings contribute; if there are
$n$ non-virtual crossings, then they divide the diagram of the knot
into $n$ components, each going from undercrossing to undercrossing,
some of which may pass through virtual crossings. In particular, the
cord algebra can be defined for virtual knots as for usual knots,
with the added stipulation that endpoints of cords can move
unobstructed through virtual crossings.

Both the cord algebra and the framed knot DGA constitute invariants
of virtual knots. In fact, they are invariants under a weaker
equivalence, that of welded knots \cite{Sat}. The fact that the cord
algebra gives an invariant of welded knots essentially follows from
the proof that it gives an invariant of usual knots; one can also
verify that the framed knot DGA is also invariant by checking the
virtual Reidemeister moves and the relevant ``forbidden'' move,
along the lines of the discussion at the end of
Section~\ref{ssec:cordDGA}. We omit the details here.

The extension of the framed knot DGA to welded knots is closely
related to an extension of the braid homomorphism $\phi$ from
Section~\ref{ssec:DGAdef} to the welded braid group \cite{FRR}. The
welded braid group $\WB_n$ adds to the usual braid group $B_n$ (with
generators $\{\sigma_k\}_{k=1}^{n-1}$) a set of generators
$\{\tau_k\}_{k=1}^{n-1}$, and relations $\tau_k^2=1$,
$\tau_k\tau_{k+1}\tau_k = \tau_{k+1}\tau_k\tau_{k+1}$,
$\tau_k\tau_{k+1}\sigma_k = \sigma_{k+1}\tau_k\tau_{k+1}$,
$\sigma_k\sigma_{k+1}\tau_k = \tau_{k+1}\sigma_k\sigma_{k+1}$, and
$\tau_k\tau_j = \tau_j\tau_k$, $\sigma_k\tau_j = \tau_j\sigma_k$ for
$|k-j|>1$. Now define $\phi_{\tau_k}:\thinspace \A_n \to \A_n$ to
act on generators $a_{ij}$ by replacing the index $k$ whenever it
appears by $k+1$ and vice versa. It is easy to check that the
relations in $\WB_n$ are satisfied and hence $\phi$ descends to a
homomorphism $\WB_n \to \Aut\A_n$. Any welded knot is the closure of
a welded braid, and one can use $\phi$ to give an alternate
definition for the framed knot DGA of a welded knot.

We remark that the fact that the framed knot DGA extends to welded
knots suggests that it may be related to knot quandles or racks.
Indeed, one presentation \cite{FR} of the fundamental rack of a
framed knot $K$ involves homotopy classes of paths in $S^3\setminus
K$ from a fixed base point in $S^3$ to a parallel copy of $K$. By
comparison, one can define the cord algebra of $K$ to be generated
by paths in $S^3\setminus K$ beginning and ending on a parallel copy
of $K$, and a cord is then a juxtaposition of two elements of the
fundamental rack. The defining relations for the cord algebra can be
expressed in terms of the binary operation on the fundamental rack,
at least if we set $\lambda=1$ in the cord algebra. A more precise
correspondence between the cord algebra and the fundamental rack
currently remains elusive, however.

\subsecdivide
%*********************************************************************
\subsection{$2$-knots in $\R^4$}
\label{ssec:2knots}

Definition~\ref{def:generalcordalg} yields a definition of the cord
algebra for embedded surfaces ($2$-knots) in $\R^4$ which extends
the cord ring from the appendix to \cite{Ng2}. We claim that this
invariant is highly nontrivial. Let $K$ be a knot in $\R^3$, and
denote by $\tilde{K}$ the spun knot (topologically $S^2$) in $\R^4$
obtained from $K$. Then $\pi_1(\R^4\setminus\tilde{K}) =
\pi_1(\R^3\setminus K)$, and $H_1(\partial\nu\tilde{K}) = \Z$ is
generated by the meridian $m$ of $K$.

It follows that the cord algebra of $\tilde{K}$ is the tensor
algebra over $\Z[\mu^{\pm 1}]$ freely generated by
$\pi_1(S^3\setminus K)$, modulo the relations (\ref{htpy1}) and
(\ref{htpy3}) from Definition~\ref{def:knotcordalg}. In other words,
its defining relations are the same as those of the cord algebra of
$K$, but without $[\gamma l]=[l\gamma]=\lambda[\gamma]$.
Equivalently, the cord algebra of $\tilde{K}$ is the cord algebra of
any long knot whose closure is $K$.

Using this formulation, one computes that the cord algebra of the
unknotted $S^2$ in $\R^4$ (the spun unknot) is the base ring
$\Z[\mu^{\pm 1}]$, while the cord algebra of the spun trefoil (with
either chirality) is $\Z[\mu^{\pm 1}][x]/((x-\mu-1)(\mu x+1))$. Some
indication of the strength of the cord algebra for spun knots comes
from the following generalization of \cite[Prop.~A.6]{Ng2}.

\begin{proposition}
The cord algebra distinguishes the spun version of any nontrivial
knot from the unknotted $S^2\subset\R^4$.
\end{proposition}

\begin{proof}
Let $K$ be a nontrivial knot. The nontriviality of the
$A$-polynomial of $K$ implies that for generic $\mu_0\in\C^*$, there
are at least two maps $HC_0(K)\otimes\C \to \C$ sending $\mu$ to
$\mu_0$ (and $\lambda$ to any nonzero complex number). It follows
that there are at least two maps from the complexified cord algebra
of $\tilde{K}$ to $\C$ sending $\mu$ to $\mu_0$. This is not the
case when $K$ is the unknot.
\end{proof}

It appears that the cord algebra may constitute an interesting
invariant for twist-spun knots as well, though computation of the
cord algebra for twist-spun knots might be a bit more involved than
for spun knots.

%*********************************************************************
%*********************************************************************

\end{document}